\documentclass[11pt]{article}

\usepackage[T1]{fontenc}
\usepackage[utf8]{inputenc}
\usepackage{lmodern}

\usepackage{amsmath,amssymb,amsthm}
\usepackage{hyperref}

\usepackage{quiver}
\usepackage{tikz}
\usepackage{tikz-cd}
\usepackage[all]{xy}
\usepackage[active]{srcltx}

\theoremstyle{plain}
\newtheorem{teo}{Theorem}[section]
\newtheorem{prop}[teo]{Proposition}
\newtheorem{lema}[teo]{Lemma}
\newtheorem{cor}[teo]{Corollary}

\theoremstyle{definition}
\newtheorem{deff}[teo]{Definition}
\newtheorem{obs}[teo]{Remark}
\newtheorem{exem}[teo]{Example}

\newtheorem*{deff*}{Definition}
\newtheorem*{teo*}{Theorem}
\newtheorem*{ack*}{Acknowledgments}
\newtheorem*{Funding*}{Funding}

\newcommand{\edg}{\color{black}}

\newcommand{\edb}{\color{black}}

\newcommand\Hom{\operatorname{Hom}}
\newcommand\GR{\operatorname{Gr}}
\newcommand\Ext{\operatorname{Ext}}
\newcommand\modd{\operatorname{mod}}

\newcommand\Fac{\operatorname{Fac}}
\newcommand\Sub{\operatorname{Sub}}
\newcommand\add{\operatorname{add}}

\newcommand{\leqt}{\preceq}
\newcommand{\geqt}{\succeq}
\newcommand{\leet}{\prec}
\newcommand{\geet}{\succ}

\numberwithin{equation}{section}

\title{Stratifying Systems via a Nested Family of Torsion Pairs}

\author{
	Edson Ribeiro Alvares\\
	Department of Mathematics\\
	Federal University of Paraná\\
	\texttt{rolo@ufpr.br}
	\and
	Matheus Vinicius dos Santos\\
	Department of Mathematics\\
	Federal University of Paraná\\
	\texttt{umchazazibalo@gmail.com}
}

\date{}

\begin{document}
	
	\maketitle
	
	\begin{center}
		\small
		This is the author's accepted manuscript of the article
		
		\vspace{0.3em}
		
		\emph{Edson Ribeiro Alvares and Matheus Vinicius dos Santos,
			``Stratifying Systems via a Nested Family of Torsion Pairs,''}
		\emph{Revista Matemática Iberoamericana}
		\textbf{42} (2026), no.~4, 1201--1226.
		
		\vspace{0.3em}
		
		The final Version of Record is available at
		
		\url{https://doi.org/10.4171/RMI/1625}
	\end{center}
	
	\vspace{1em}
	
	\begin{center}
		
		\vspace{0.5em}
		
	\end{center}
	
	\begin{abstract}
		In this paper{\edg,} we introduce the notion of a nested family of torsion pairs. We show that every stratifying system induces a nested family of torsion pairs and, conversely, that such a family together with a suitable object gives rise to a stratifying system; in particular, every stratifying system arises in this way. As an application, we present a stratifying system of infinite size that cannot be indexed by $(\mathbb{N}, \leq)$, where $\leq$ is the natural order.
	\end{abstract}

\section{Introduction}

{\edg Stratifying systems were introduced by K. Erdmann and C. Sáenz in \cite{definicaosistemaestratificante} as a tool to construct, for a finite-dimensional algebra $A$ over an algebraically closed field $k$, an $A$-module whose endomorphism ring is standardly stratified.} This original definition is known in the literature as an Ext-injective stratifying system. {\edg Subsequently, E. Marcos, O. Mendoza, and C. Sáenz proposed in \cite{stratifyingsimple} a new definition. In the present work, we consider a generalization where the indexing set is an arbitrary totally ordered set $(G, \leq)$ (see Definition \ref{definicao}). The definition of Marcos, Mendoza and Sáenz is equivalent to the one presented here in the special case where $G$ is the finite set $\{1, \dots, t\}$.}

While the category of modules
filtered by an exceptional sequence is equivalent to the subcategory of modules filtered by
standard modules in the category of a quasi-hereditary algebra \cite[Theorem 1.6, 1.7]{definicaosistemaestratificante}, the
category of modules filtered by a stratifying system is likewise equivalent to the subcategory of modules filtered by standard modules in the category of modules over a standardly
stratified algebra \cite[Theorem 1.6]{definicaosistemaestratificante}.

{\edg Furthermore, the homological properties of these stratifying systems play an important role in the context of quasi-hereditary and standardly stratified algebras \cite{Cline1988, Quasi}. Consequently, stratifying systems have attracted considerable interest from researchers in representation theory.} Given a stratifying system, it is well known that the full subcategory consisting of the zero module and all modules admitting a filtration by this system is closed under direct sums and direct summands. In this regard, the pioneering work of V. Dlab and C. Ringel \cite{MR1211481} serves as a foundational reference. Another important approach was presented by I. Ágoston, V. Dlab and E. Lukács in \cite{AGOSTON20084177}, where they proved that, given a finite-dimensional algebra, there exists a unique basic algebra which is filtered by standard modules, and the subcategories of filtered modules for both algebras are equivalent. 

In the structural analysis of the subcategory of filtered modules, the existence of a stratifying system plays a central role. Ringel \cite{MR1128706} showed that this subcategory is functorially finite in the sense of Auslander--Smalø \cite{MR1097029} and that it admits almost split sequences. More recently, Bautista, P\'erez, and Salmer\'on \cite{MR4736629} furthered this line of research by describing the structure of such sequences and proving that the category of filtered modules can only have one of the following representation types: finite, tame, or wild.

Given a set of modules, to verify whether it forms a stratifying system, it is often sufficient to directly check whether certain morphism properties and extension conditions hold. Nevertheless, an important question is how to find a stratifying system. Despite significant advances, with the exception of hereditary algebras \cite{Cadavid1, Cadavid2}, little is known about the existence of non-trivial stratifying systems (defined by the standard modules) in arbitrary algebras. An important first step in this direction was made by O. Mendoza and H. Treffinger \cite{stratifying}, who demonstrated that it is possible to construct new stratifying systems with {\edg indecomposable direct summand}s from a certain quotient of a $\tau$-rigid module \cite[Theorem 3.4]{stratifying}. However, since a basic $\tau$-rigid module $M$ has at most $\vert A \vert$ direct summands \cite[Proposition 1.3]{reiten}, the size of a stratifying system induced by $M$ is bounded by $\vert A \vert$ \cite[Theorem 3.4]{stratifying}. This technique forms the foundation for the central idea we introduce in this paper. Specifically, from these stratifying systems, we can generate what we call \textit{nested families of torsion pairs}, with the stratifying system being intricately linked to this family in a special way.

\subsection*{\bfseries Overview of the text structure and main results}

In the main result of this paper{\edg,} we show that the construction of stratifying systems can be carried out in a more general setting by considering the concept of a nested family of torsion pairs. This is a collection of torsion pairs $\Gamma = \{ ( \mathcal{T}_{k}, \mathcal{F}_{k}) \}_{k \in G}$ in an abelian category $\mathcal{A}$, indexed by a totally ordered set $(G, \leq)$, satisfying $\mathcal{T}_{k} \supsetneq \mathcal{T}_{l}$ whenever $k < l$, or equivalently, $\mathcal{F}_{k} \subsetneq \mathcal{F}_{l}$ whenever $k < l$ (Definition \ref{nested}).

To develop this construction, we work within the category \(\GR_{G}(\mathcal{A})\) of \(G\)-graded objects in \(\mathcal{A}\), which provides a natural framework for introducing the key subcategories \(\mathcal{M}(\Gamma)\) and \(\mathcal{N}(\Gamma)\) (Definition \ref{compatibility}). The subcategory \(\mathcal{M}(\Gamma)\) consists of objects that lie in specific torsion parts of the nested family and satisfy certain compatibility conditions across indices; similarly, \(\mathcal{N}(\Gamma)\) is defined using the torsion-free parts. These subcategories are fundamental to constructing pre-stratifying systems, a concept that generalizes stratifying systems {\edg where the indexed objects need not be indecomposable, yet essential orthogonality conditions are maintained} (Definition \ref{definicao}). In particular, by choosing an {\edg indecomposable direct summand} of each object in a pre-stratifying system, one obtains a stratifying system.

The primary motivation for introducing pre-stratifying systems is to provide a unified organizational framework. For each nested family of torsion pairs $\Gamma$, we identify specific objects in $\mathcal{M}(\Gamma)$ (respectively, $\mathcal{N}(\Gamma)$) and associate to them objects in the more structured subcategories $\mathcal{M}^{*}(\Gamma)$ (respectively $\mathcal{N}^{*}(\Gamma)$). This association is realized via the \textit{stratum} and \textit{substratum} functors, $\mathfrak{F} \colon \mathcal{M}(\Gamma) \to \mathcal{M}^{*}(\Gamma)$ and $\mathfrak{T} \colon \mathcal{N}(\Gamma) \to \mathcal{N}^{*}(\Gamma)$, whose definitions depend explicitly on the choice of $\Gamma$ (Definition \ref{stratum}).

With this established, the underlying ideas of the text can be summarized as follows. In one of the main results, we prove that objects in the category $\GR_G(\mathcal{A})$ satisfying a certain orthogonality property induce two distinct types of nested families of torsion pairs (Theorem \ref{galo}). Consequently, given a pre-stratifying system, we can associate to it an object in $\GR_G(\mathcal{A})$ which, by this result, induces two distinct nested families of torsion pairs (Corollary \ref{nestedfamilyinducedbystrat}).

Conversely, from a nested family of torsion pairs, we organize objects into a subcategory of \(\mathcal{M}(\Gamma)\) (respectively \(\mathcal{N}(\Gamma)\)) denoted \(\mathcal{M}^{\dagger}(\Gamma)\) (respectively \(\mathcal{N}^{\dagger}(\Gamma)\)), whose objects satisfy a specific property (Definition \ref{aqui}). From each object in these categories, we can induce a pre-stratifying system (Theorem \ref{teoresma1}).

These are the central ideas of this paper. With them, we show that the result of Treffinger and Mendoza \cite[Theorem 3.4]{stratifying} forms part of a more general context in which one can induce pre-stratifying systems (see Corollary \ref{teoremadoscara}).

It is also important to highlight another aspect that motivated us to present the more general definition of stratifying systems (see Definition \ref{definicao}). It is known that even in representation-finite algebras it is possible to find stratifying systems of size greater than $\vert A \vert$ \cite[Remark 2.7]{stratifyingsimple}.  In \cite{sistemainfinito}, H. Treffinger presents a stratifying system indexed by $\mathbb{N}$, demonstrating that a stratifying system can indeed have infinite size. Since there does not always exist an order-preserving bijection between totally ordered sets $(G, \leq)$ and $(H, \preceq)$ (even if $G$ and $H$ have the same cardinality), the possibility of stratifying systems indexed by infinite sets other than $\mathbb{N}$ further justifies the need for this more general definition.

\subsection*{\bfseries Outline of the Paper}

We begin in Section \ref{sec2} by fixing notation and recalling necessary background on torsion pairs and $\tau$-tilting theory.

We present our general framework in Section \ref{sec3}, introducing nested families of torsion pairs and the categories $\mathcal{M}(\Gamma)$ and $\mathcal{N}(\Gamma)$ whose objects reflect the nested structure of $\Gamma$. The stratum and substratum functors are then defined, mapping objects from $\mathcal{M}(\Gamma)$ to $\mathcal{M}^{*}(\Gamma)$ and from $\mathcal{N}(\Gamma)$ to $\mathcal{N}^{*}(\Gamma)$, respectively (Lemma \ref{lesma1}).

In Section \ref{sec4}, we establish the {\edg main} connections. We define pre-stratifying systems and show how they interact with our torsion framework. Theorem \ref{galo} demonstrates that objects in $\GR_G(\mathcal{A})$ with an orthogonality property naturally induce two distinct nested families. Corollary \ref{nestedfamilyinducedbystrat} ensures that every pre-stratifying system provides such an object. In the opposite direction, we define the subcategories $\mathcal{M}^{\dagger}(\Gamma)$ and $\mathcal{N}^{\dagger}(\Gamma)$ whose objects are precisely those that induce pre-stratifying systems (Theorem \ref{teoresma1}). Finally, Theorem \ref{cubo} completes this bridge by showing that every pre-stratifying system arises in this way.

The subsequent sections explore consequences and applications. In Section \ref{sec5}, we show that the construction of stratifying systems from $\tau$-rigid modules \cite{stratifying} is a particular case of our approach (Corollary \ref{teoremadoscara}). Section \ref{sec6} introduces an order relation among torsion families, addressing the case where an object belongs to $\mathcal{M}^{\dagger}(\Gamma)$ for multiple families $\Gamma$. Section \ref{sec7} exhibits a stratifying system that cannot be indexed by the natural numbers with the natural order, validating the generality of our definition indexed by an arbitrary ordered set $G$. Finally, Section \ref{x1x1} provides a detailed example illustrating our main constructions.

\section{Setting and Background}\label{sec2}

Throughout this paper, $\mathcal{A}$ denotes a $k$-linear, Hom-finite, skeletally small abelian category, and $A$ denotes a finite-dimensional $k$-algebra over an algebraically closed field $k$.
We denote by $\modd A$ the category of finitely generated right $A$-modules. For $M \in \modd A$, the number of pairwise non-isomorphic {\edg indecomposable direct summand}s of $M$ is denoted by $\vert M \vert$. If $M$ does not have two isomorphic direct summands, we say that $M$ is a basic module. We write $X \mid M$ to indicate that $X$ is a direct summand of $M$.
The Auslander-Reiten translation is denoted by $\tau$ and $D = \Hom_k(-,k)$ is the duality functor. $G$ will denote a set with a total order $\leq$, and $G^{op}$ will denote the same set with the total order $\leq^{op}$, where $x \leq^{op} y$ in $G^{op}$ if and only if $y \leq x$ in $G$. $\mathbb{I}_n$ will denote the set $\{1, \cdots, n\}$ with the natural order.

For a given class $\mathcal{C} \subseteq \mathcal{A}$, we define 
\begin{align*}
	\mathcal{C}^{\bot} := \{M \in \mathcal{A}\ \mid \  \Hom_{\mathcal{A}}(-,M)\vert_{\mathcal{C}}= 0\}.\\
	{}^{\bot}\mathcal{C} := \{M \in \mathcal{A}\ \mid \ \Hom_{\mathcal{A}}(M,-)\vert_{\mathcal{C}} = 0\}.
\end{align*}
If $\mathcal{C} = \{M\}$ for some $M \in \mathcal{A}$, we will denote $\{M\}^{\bot}$ and ${}^{\bot}\{M\}$ simply by $M^{\bot}$ and ${}^{\bot}M$, respectively. An object $M \in \mathcal{C}$ is called {\bfseries Ext-projective} in $\mathcal{C}$ if $\Ext_{\mathcal{A}}^1(M,-)\vert_{\mathcal{C}}=0$. Dually, it is called {\bfseries Ext-injective} in $\mathcal{C}$ if $\Ext_{\mathcal{A}}^1(-, M)\vert_{\mathcal{C}} = 0$.

For a given $M \in \modd A$ we define
\begin{align*}
	\Fac(M) := \{X \in \modd A \ \mid \ \exists \hbox{ an epimorphism } M^{\oplus n} \rightarrow X \hbox{ for some } n \in \mathbb{N}\}.\\
	\Sub(M) := \{X \in \modd A \ \mid \ \exists \hbox{ a monomorphism } X \rightarrow M^{\oplus n}\hbox{ for some } n \in \mathbb{N}\}.
\end{align*}

\subsection*{\bf Torsion Pairs} The results presented in this paragraph can be found in \cite[Chapter VI]{livroazul} and in \cite{Bruce}.

A pair $(\mathcal{T}, \mathcal{F})$ of full subcategories of an abelian category $\mathcal{A}$ is called a {\bfseries torsion pair in $\boldsymbol{\mathcal{A}}$} if $\Hom_{\mathcal{A}}(T, F) = 0$ for all $T \in \mathcal{T}$ and $F \in \mathcal{F}$, and if for every object $M \in \mathcal{A}$ there exists a short exact sequence
\[
0 \longrightarrow t(M) \longrightarrow M \longrightarrow f(M) \longrightarrow 0,
\]
with $t(M) \in \mathcal{T}$ and $f(M) \in \mathcal{F}$. 
The subcategories $\mathcal{T}$ and $\mathcal{F}$ are called the {\bf torsion class} and {\bf torsion-free class}, respectively. 
Moreover, this exact sequence gives rise to additive functors $t, f: \mathcal{A} \longrightarrow \mathcal{A}$, called the {\bf torsion functor} and the {\bf torsion-free functor}, respectively. We have that $t(M) \longrightarrow M$ is a $\mathcal{T}$-right approximation of $M$, and $M \longrightarrow f(M)$ is a $\mathcal{F}$-left approximation of $M$. In particular, $\mathcal{T}$ is always contravariantly finite and $\mathcal{F}$ is always covariantly finite. Moreover, since $t(M) \to M$ is a monomorphism and $M \to f(M)$ is an epimorphism, these approximations are minimal.

The torsion class $\mathcal{T}$ is closed under quotients and extensions, while the torsion-free class $\mathcal{F}$ is closed under subobjects and extensions.
A torsion pair $(\mathcal{T}, \mathcal{F})$ is called {\bf splitting} if, for every $M \in \mathcal{A}$, the corresponding short exact sequence splits.

It is known that if \(\mathcal{T} = \Fac(M)\) is a torsion class for some \(M \in \modd A\), then \(\mathcal{F} = M^{\bot}\); similarly, if \(\mathcal{F} = \Sub(M)\) is a torsion-free class, then \(\mathcal{T} = {}^{\bot}M\). 

If $\mathcal{C}$ is a class in an abelian length category $\mathcal{A}$, then the smallest torsion class containing $\mathcal{C}$ is $\mathsf{T}(\mathcal{C}) := {}^{\bot}(\mathcal{C}^{\bot})$. This yields the torsion pair $(\mathsf{T}(\mathcal{C}), \mathcal{C}^{\bot})$. Analogously, the smallest torsion-free class containing $\mathcal{C}$ is $\mathsf{F}(\mathcal{C}) := ({}^{\bot}\mathcal{C})^{\bot}$, which yields the torsion pair $({}^{\bot}\mathcal{C}, \mathsf{F}(\mathcal{C}))$. {\edb In particular}, if $\mathcal{T}$ is a torsion class containing $\mathcal{C}$, then $\mathsf{T}(\mathcal{C}) \subseteq \mathcal{T}$. Similarly, if $\mathcal{F}$ is a torsion-free class containing $\mathcal{C}$, then $\mathsf{F}(\mathcal{C}) \subseteq \mathcal{F}$.

\subsection*{\bf $\tau$-Tilting Theory}

We remind the reader of the definition of a $\tau$-rigid module.
An $A$-module $M$ is said to be {\bf $ \boldsymbol\tau$-rigid} if $\Hom_A(M, \tau M) = 0$; similarly, it is said to be {\bf $ \boldsymbol\tau^{-}$-rigid} if $\Hom_A(\tau^{-1}M, M) = 0$. If $M$ is a $\tau$-rigid {\edb (or $\tau^{-}$-rigid)} module, then $\vert M \vert \leq \vert A_A \vert$ \cite[Proposition 1.3]{reiten}.

We will make use of the following results from $\tau$-tilting theory throughout the article.

\begin{prop}\cite[Proposition 3.2]{stratifying}\label{TFantigo} Let $M$ be a basic non-zero $\tau$-rigid module. Then there exists a decomposition $M = \bigoplus\limits_{i=1}^{t} M_i$ into the direct sum of indecomposable $A$-modules such that $M_i \notin \Fac(\bigoplus\limits_{j>i}M_j)$.
\end{prop}

Such a decomposition is called {\bf torsion-free admissible}.

\begin{prop}\cite[Proposition 5.8]{almostsplit}\label{rigidoextproj}
	The following statements are equivalent for a pair of modules $M, N \in \modd A$.    
	\begin{enumerate}
		\item $\Ext_A^1(M, N'') = 0$ if $N'' \in \Fac(N)$.
		\item $\Hom_A(N, \tau M) = 0$.
	\end{enumerate}
\end{prop}
In particular, if $M$ is $\tau$-rigid, then $M$ is Ext-projective in $\Fac(M)$.

\begin{prop}\cite[Theorem 5.10]{almostsplit} \label{taurigidoalmostsplit}
	If $M$ is $\tau$-rigid, then $\Fac(M)$ is a functorially finite torsion class.    
\end{prop}

\begin{lema}\cite[Lemma 4.6]{exceptional}\label{ind}
	If $X$ is an indecomposable $A$-module such that $X\bigoplus U$ is $\tau$-rigid, then either $f(X)$ is indecomposable or $f(X) = 0$, where $f$ is the torsion-free functor with respect to the torsion pair $(\Fac(U), U^{\bot})$. {\edb Moreover}, we have $f(X) = 0$ if and only if $X$ is in $\Fac(U)$.    
\end{lema}

\section{Nested Families of Torsion Pairs}\label{sec3}

In this section, we introduce the concept of nested families of torsion pairs, which will play a central role in the study of stratifying systems in Section 4. 

\subsection*{\bf Nested Families of Torsion Pairs}
\begin{deff}[Nested family of torsion pairs]
	\label{nested}
	Let $\mathcal{A}$ be an abelian category. A {\bf nested family of torsion pairs in $\boldsymbol{\mathcal{A}}$} (nested family, for short) is a set $\Gamma = \{(\mathcal{T}_k, \mathcal{F}_k) \}_{k \in G}$ of torsion pairs in $\mathcal{A}$  indexed by $(G, \leq)$ satisfying $\mathcal{T}_k \supsetneq \mathcal{T}_{l}$ if $k < l$, or equivalently, $\mathcal{F}_k \subsetneq \mathcal{F}_{l}$ if $k < l$, with $k, l \in G$.
\end{deff}

\begin{obs}
	We will denote by $t_k$ and $f_k$ the torsion and torsion-free functors with respect to the torsion pair $(\mathcal{T}_k, \mathcal{F}_k)$, respectively.
\end{obs}

Let $\Gamma = \{(\mathcal{T}_k, \mathcal{F}_k)\}_{k \in G}$ be a nested family of torsion pairs in $\mathcal{A}$. An easy calculation shows that if $(\mathcal{T}_k, \mathcal{F}_k)$ is a torsion pair in $\mathcal{A}$, then $(\mathcal{F}_k^{op}, \mathcal{T}_k^{op})$ is a torsion pair in $\mathcal{A}^{op}$, and consequently $\Gamma^{op} := \{(\mathcal{F}_k^{op}, \mathcal{T}_k^{op})\}_{k \in G^{op}}$ is a nested family in $\mathcal{A}^{op}$.

Let $\mathcal{A}$ be an abelian category and $G$ a set. We denote by $\GR_G(\mathcal{A})$ the additive category of $G$-graded objects over $\mathcal{A}$, whose objects $M_{\bullet} = (M_j)_{j \in G}\in \GR_G(\mathcal{A})$ are sequences of objects in $\mathcal{A}$ indexed by $G$, that is, $M_j \in \mathcal{A}$ for all $j \in G$. 

For $M_{\bullet} = (M_j)_{j \in G}$ and $N_{\bullet} = (N_j)_{j \in G}$ in $\GR_G(\mathcal{A})$, a morphism $f : M_{\bullet} \longrightarrow N_{\bullet}$ is a collection $f = (f_j)_{j \in G}$ such that $f_j : M_j \longrightarrow N_j$ is a morphism in $\mathcal{A}$ for every $j \in G$. The composition of morphisms in $\GR_G(\mathcal{A})$ is defined componentwise, that is, if $g: N_{\bullet} \longrightarrow P_{\bullet}$ is a morphism such that $g = (g_j)_{j\in G}$, then $(g \circ f) := (g_j \circ f_j)_{j\in G}$. Note that $\GR_G(\mathcal{A})$ is an additive category with zero object $0_{\bullet} = (0)_{j \in G}$. Moreover, if $G$ is finite, then $\GR_G(\mathcal{A})$ can be identified with a (not necessarily full) subcategory of $\mathcal{A}$.

Given a nested family of torsion pairs $\Gamma$ in $\mathcal{A}$ indexed by $G$, the objects $M_{\bullet}$ in $\GR_G(\mathcal{A})$ that exhibit a certain compatibility with the torsion (or torsion-free) classes of $\Gamma$ will play an important role.

\begin{deff}\label{compatibility}
	Let $\Gamma = \{ (\mathcal{T}_k, \mathcal{F}_k) \}_{k \in G}$ be a nested family of torsion pairs in $\mathcal{A}$. We define $\mathcal{M}(\Gamma)$ as the full subcategory of $\GR_G(\mathcal{A})$ consisting of $0_{\bullet}$ and the objects $M_{\bullet} = (M_k)_{k \in G}$ satisfying the following conditions for all $k \in G$:
	\begin{enumerate}
		\item[t1.] $M_k \in \mathcal{A}$ is non-zero;
		\item[t2.] $M_k \in \mathcal{T}_k$;
		\item[t3.] $M_k \notin \mathcal{T}_j$ whenever $j > k$.
	\end{enumerate}
	Furthermore, we define $\mathcal{M}^{*}(\Gamma)$ as the full subcategory of $\mathcal{M}(\Gamma)$ consisting of $0_{\bullet}$ and the objects $M_{\bullet} = (M_k)_{k \in G}\in\mathcal{M}(\Gamma)$ that satisfy the condition that $M_k \in \mathcal{F}_j$ whenever $j > k$, for all $k \in G$.	
\end{deff}

{\edb The dual construction yields the subcategories $\mathcal{N}(\Gamma)$ and $\mathcal{N}^*(\Gamma)$, defined as follows.}
\begin{deff}
	
	Let $\Gamma = \{ (\mathcal{T}_k, \mathcal{F}_k) \}_{k \in G}$ be a nested family of torsion pairs in $\mathcal{A}$. We define $\mathcal{N}(\Gamma)$ as the full subcategory of $\GR_G(\mathcal{A})$ consisting of $0_{\bullet}$ and the objects $N_{\bullet} = (N_k)_{k \in G}$ satisfying the following conditions for all $k \in G$:
	\begin{enumerate}
		\item[f1.] $N_k \in \mathcal{A}$ is non-zero;
		\item[f2.] $N_k \in \mathcal{F}_k$;
		\item[f3.] $N_k \notin \mathcal{F}_j$ whenever $j < k$.
	\end{enumerate}
	Additionally, we define $\mathcal{N}^{*}(\Gamma)$ as the full subcategory of $\mathcal{N}(\Gamma)$ consisting of $0_{\bullet}$ and the objects $N_{\bullet} = (N_k)_{k \in G} \in \mathcal{N}(\Gamma)$ that satisfy the condition that $N_k \in \mathcal{T}_j$ whenever $j < k$, for all $k \in G$.
	
\end{deff}

We have that the category $\mathcal{M}(\Gamma)$ is an additive subcategory of $\GR_G(\mathcal{A})$. Indeed, if $M_{\bullet} = (M_k)_{k \in G}$ and $N_{\bullet} = (N_k)_{k \in G}$ are objects in $\mathcal{M}(\Gamma)$, then for every $k \in G$ we have that $M_k \bigoplus N_k$ is non-zero, and $M_k \bigoplus N_k \in \mathcal{T}_k$, since $\mathcal{T}_k$ is closed under direct sums. Moreover, \( M_k \bigoplus N_k \notin \mathcal{T}_j \) for \( j > k \), since neither \( M_k \) nor \( N_k \) lies in \( \mathcal{T}_j \), and \( \mathcal{T}_j \) is closed under direct summands. Consequently, $M_{\bullet} \bigoplus N_{\bullet} = (M_k \bigoplus N_k)_{k \in G} \in \mathcal{M}(\Gamma)$. Analogously, one proves that $\mathcal{M}^{*}(\Gamma)$, $\mathcal{N}(\Gamma)$, and $\mathcal{N}^{*}(\Gamma)$ are additive subcategories of $\GR_G(\mathcal{A})$.

Let $\Gamma = \{(\mathcal{T}_k, \mathcal{F}_k)\}_{k \in G}$ be a nested family in the abelian category $\mathcal{A}$. The following pairs of categories form contravariant equivalences (or dualities):
\[
(\mathcal{M}(\Gamma), \mathcal{N}(\Gamma^{op})), (\mathcal{N}(\Gamma), \mathcal{M}(\Gamma^{op})), (\mathcal{M}^{*}(\Gamma), \mathcal{N}^{*}(\Gamma^{op})), (\mathcal{N}^{*}(\Gamma), \mathcal{M}^{*}(\Gamma^{op})).
\]

Indeed, let $M_{\bullet} =(M_k)_{k \in G} \in\mathcal{M}(\Gamma)$. We have that $M_k \notin \mathcal{T}_j$ if $j > k$ for all $k \in G$, which is equivalent to $M_k \notin \mathcal{T}_j$ if $j <^{op} k$ for all $k \in G$. This shows that $M_{\bullet} \in\mathcal{N}(\Gamma^{op})$. (Note that $\mathcal{T}_j^{op}$ is the torsion-free part of the pair $(\mathcal{F}_j^{op}, \mathcal{T}_j^{op})$ in $\mathcal{A}^{op}$). Moreover, we have that $\mathcal{M}(\Gamma) \subseteq \mathcal{N}(\Gamma^{op})$; similarly, $\mathcal{N}(\Gamma^{op}) \subseteq \mathcal{M}(\Gamma)$ can be shown. The other cases follow analogously.

The categories $\mathcal{M}(\Gamma)$ and $\mathcal{N}(\Gamma)$ will play a dual role in the development of the paper. For the sake of completeness, we will state both versions of each result; however, we will only prove the ones related to $\mathcal{M}(\Gamma)$, as the other case is analogous.

\begin{exem}\label{exem512512}
	We will show that for any totally ordered set $(G, \leq)$, there exists a nested family indexed by $G$ in $\modd A$ for some finite-dimensional algebra $A$ over an algebraically closed field $k$.
	
	{\edg Let \( G \) be an infinite totally ordered set}. If \( \hbox{card}\ G = \hbox{card}\ \mathbb{N} \), take \( k = \overline{\mathbb{Q}} \), where \(\overline{\mathbb{Q}} \subseteq \mathbb{C}\) denotes the field of algebraic numbers, that is, the set of roots of non-zero polynomials with integer coefficients; if \( \hbox{card}\ G > \hbox{card}\ \mathbb{N} \), consider the integral domain \( \overline{\mathbb{Q}}[G] \) defined as the polynomial ring over \( \overline{\mathbb{Q}} \) in the set of commuting indeterminates indexed by \( G \). We have \( \hbox{card}\ \overline{\mathbb{Q}}[G] = \hbox{card}\ G \). We construct the field of fractions \( \hbox{frac}\ \overline{\mathbb{Q}}[G] \) and take its algebraic closure \( k \). Both constructions preserve the cardinality of \( G \), so that \( \hbox{card}\ k = \hbox{card}\ G \).  
	
	Since \( G \) and \( k \) have the same cardinality, there exists a bijection \( \varphi: G \longrightarrow k \). Through this bijection, we can transfer the total order of \( G \) to \( k \), defining \( \lambda_1 < \lambda_2 \) in \( k \) if and only if \( \varphi^{-1}(\lambda_1) < \varphi^{-1}(\lambda_2) \) in \( G \).

	Let $ A = kQ$ be the Kronecker algebra. For every $\lambda \in k$, consider the module $ M_{\lambda} $ whose representation is given by
	\[
	\begin{tikzcd}
		k & k
		\arrow["\lambda"', shift right, from=1-2, to=1-1]
		\arrow["1", shift left, from=1-2, to=1-1]
	\end{tikzcd}
	\]
	and consider for every $i \in G$ the torsion pairs $(\mathcal{T}_i, \mathcal{F}_i)$, where 
\[
\mathcal{T}_i = \mathsf{T}(\{M_{\lambda} \mid \lambda \geq \varphi(i)\}) 
\quad\text{and}\quad 
\mathcal{F}_i = (\{M_{\lambda} \mid \lambda \geq \varphi(i)\})^{\bot}.\]

	We have that $\Gamma = \{(\mathcal{T}_i, \mathcal{F}_i)\}_{i \in G}$ is a nested family. Indeed, it is immediate from the definition that if $j < l${\edg,} then $\mathcal{T}_{j} \supseteq \mathcal{T}_{l}$. To show that the inclusion is proper, we just notice that $M_{\lambda} \in \mathcal{T}_{\varphi^{-1}(\lambda)}$. However, since $\Hom_A(M_{\rho}, M_{\lambda}) = 0 $ if $\lambda \neq \rho$, we have that $M_{\lambda} \in \mathcal{F}_{\varphi^{-1}(\rho)}$ implies that $M_{\lambda} \notin \mathcal{T}_{\varphi^{-1}(\rho)}$ whenever $\rho > \lambda$. Therefore, $\mathcal{T}_{j} \supsetneq \mathcal{T}_{l}$.
	
	If $\hbox{card}\ G = n$, we can follow the same construction over a finite subset of $\overline{\mathbb{Q}}$.
\end{exem}

\subsection*{\bf Stratum and Substratum}

Let $\Gamma = \{(\mathcal{T}_k, \mathcal{F}_k)\}_{k \in G}$ be a nested family of torsion pairs in an abelian length category $\mathcal{A}$. We construct a functor  
\[  
\mathfrak{F}: \mathcal{M}(\Gamma) \longrightarrow \mathcal{M}^{*}(\Gamma),  
\]  
such that for any $M_{\bullet} = (M_k)_{k \in G} \in \mathcal{M}(\Gamma)$ we have $\mathfrak{F}(M_{\bullet}) = (\mathfrak{F}_k(M_{\bullet}))_{k \in G}$,  
{\edb where each $\mathfrak{F}_k(M_{\bullet})$ is the torsion-free part of $M_k$ with respect to the torsion pair $({}^{\bot}(\bigcap\limits_{j>k} \mathcal{F}_j),\bigcap\limits_{j>k} \mathcal{F}_j)$.}

To this end, we adopt the following convention.
Fix a skeleton $\mathcal{S}$ of $\mathcal{A}$, that is, a full subcategory of $\mathcal{A}$ such that for every object $M \in \mathcal{A}$ there exists a unique object $M' \in \mathcal{S}$ with $M \cong M'$.
If $\Gamma = \{(\mathcal{T}_k, \mathcal{F}_k)\}_{k \in G}$ is a nested family of torsion pairs, we assume, without loss of generality, that for every object $M \in \mathcal{A}$ and every $k \in G$, both $t_k(M)$ and $f_k(M)$ lie in $\mathcal{S}$. Moreover, we assume, again without loss of generality, that if for some object $M \in \mathcal{A}$ we have $t_k(M) = t_i(M)$ and $f_k(M) = f_i(M)$ for some $i, k \in G$, then the corresponding canonical exact sequences coincide.
This can be done since the morphisms $i_j:t_j(M) \longrightarrow M$ and $\pi_j: M \longrightarrow f_j(M)$, for $j \in \{i, k\}$, are respectively the minimal {\edg right} $\mathcal{T}_j$-approximation and the minimal {\edg left} $\mathcal{F}_j$-approximation of $M$.

{\edb We begin by defining functors $\mathfrak{F}_k : \mathcal{M}(\Gamma) \longrightarrow \mathcal{A}$. Let $\Gamma = \{(\mathcal{T}_k, \mathcal{F}_k)\}_{k \in G}$ be a nested family of torsion pairs. For every $k \in G$, we have that
	\[
	\left({}^{\bot}\!\left(\bigcap\limits_{j>k} \mathcal{F}_j\right), \bigcap\limits_{j>k} \mathcal{F}_j\right)
	\]
	is a torsion pair in $\mathcal{A}$ \cite[Corollary 3.5, Corollary 3.6]{Bruce}. Here, we set $\bigcap\limits_{j>k} \mathcal{F}_j = \mathcal{A}$ if $k$ is the maximum element of $G$ (if such an element exists).
	
	For $M_{\bullet} = (M_j)_{j \in G} \in \mathcal{M}(\Gamma)$, we define $\mathfrak{F}_k(M_{\bullet})$ as the torsion-free part of $M_k$ with respect to the torsion pair $({}^{\bot}(\bigcap\limits_{j>k} \mathcal{F}_j),\bigcap\limits_{j>k} \mathcal{F}_j)$.
	
	This yields the following explicit description:
	
	\[
	\mathfrak{F}_k(M_{\bullet}) :=
	\begin{cases}
		f_m(M_k), & \text{if $k$ is not the maximum element of $G$,} \\
		M_k, & \text{if $k$ is the maximum element of $G$.}
	\end{cases}
	\]
	where $m$ is any index such that the set $\{t_j(M_k) \mid k < j \leq m\}$ has exactly one element.}

{\edb Indeed, if $k$ is the maximum element of $G$, then $\bigcap\limits_{j>k} \mathcal{F}_j = \mathcal{A}$, and therefore $\mathfrak{F}_k(M_{\bullet}) = M_k$. Assume henceforth that $k$ is not the maximum element of $G$ and consider the set
\[
\{t_j(M_k) \mid j > k\},
\]
which consists of subobjects of $M_k$.} {\edg Note that $t_i(M_k)$ is a subobject of $t_j(M_k)$ whenever $i>j$. Indeed, there exist monomorphisms $h_i \colon t_i(M_k) \longrightarrow M_k$ and $h_j \colon t_j(M_k) \longrightarrow M_k$, where $h_i$ and $h_j$ are minimal right $\mathcal{T}_i$-approximation and minimal right $\mathcal{T}_j$-approximation of $M_k$, respectively. Since $\mathcal{T}_i \subsetneq \mathcal{T}_j$, it follows that $t_i(M_k) \in \mathcal{T}_j$, and hence there exists a morphism
\[
g \colon t_i(M_k) \longrightarrow t_j(M_k)
\]
such that $h_i = h_j \circ g$. Consequently, $g$ is a monomorphism.}

Since $\mathcal{A}$ is an abelian length category, every object $M_k$ is both artinian and noetherian. Hence, the set $\{t_j(M_k) \mid j > k\}$ is finite (under our convention that if $t_i(M_k) \cong t_j(M_k)$, then $t_i(M_k) = t_j(M_k)$). Therefore, there exists at least an index $m > k$ such that the set $\{t_j(M_k) \mid k < j \leq m\}$ consists of a single element (note that if an index $m$ satisfies this property, then every $n$ with $k < n < m$ also satisfies it, and a maximal such $m$ need not exist). Consequently, the set $\{f_j(M_k) \mid k < j \leq m\}$ also consists of a single element.

{\edb To see that $\mathfrak{F}_k(M_{\bullet}) = f_m(M_k)$, that is, that the torsion-free part of $M_k$ with respect to the torsion pair $({}^{\bot}(\bigcap\limits_{j>k} \mathcal{F}_j),\bigcap\limits_{j>k} \mathcal{F}_j)$ is given by $f_m(M_k)$, it suffices to show that $f_m(M_k) \in \bigcap\limits_{j>k} \mathcal{F}_j$ and $t_m(M_k) \in {}^{\bot}(\bigcap\limits_{j>k} \mathcal{F}_j)$, since the canonical exact sequence of $M_k$ with respect to the torsion pair $$({}^{\bot}(\bigcap\limits_{j>k} \mathcal{F}_j),\bigcap\limits_{j>k} \mathcal{F}_j)$$ is unique.
	
	Note that $f_m(M_k) \in \mathcal{F}_j$ for all $j>k$, since
	$$
	\{f_m(M_k)\} = \{f_j(M_k) \mid k < j \leq m\}
	\Longrightarrow f_m(M_k) \in \mathcal{F}_j \text{ for } k < j \leq m,
	$$
	and hence $f_m(M_k) \in \mathcal{F}_j$ for all $j>k$, as $\mathcal{F}_i \subsetneq \mathcal{F}_j$ whenever $i<j$.
	Therefore,	
	$$f_m(M_k) \in \bigcap_{j>k} \mathcal{F}_j.$$

	Similarly, since $t_m(M_k) \in \mathcal{T}_j$ for all $k < j \leq m$, it follows that $$\Hom_{\mathcal{A}}(t_m(M_k),-)\vert_{\mathcal{F}_j} = 0$$ for all $k < j \leq m$, and consequently
	$$
	t_m(M_k) \in {}^{\bot}\left(\bigcap\limits_{k < j \leq m} \mathcal{F}_j\right) = {}^{\bot}\left(\bigcap\limits_{j>k} \mathcal{F}_j\right),
	$$
	since $\bigcap\limits_{k < j \leq m} \mathcal{F}_j = \bigcap\limits_{j>k} \mathcal{F}_j$, as $\mathcal{F}_i \subsetneq \mathcal{F}_j$ whenever $i<j$.

\begin{obs} 
	\begin{enumerate}
		\item Since $\mathcal{F}_i \subsetneq \mathcal{F}_j$ whenever $i<j$, one has
		\[
		{}^{\bot}\!\left(\bigcap\limits_{j>k} \mathcal{F}_j\right) = \bigcup\limits_{j>k} \mathcal{T}_j.
		\]
		Here, we set $\bigcup\limits_{j>k} \mathcal{T}_j = \add{0}$ if $k$ is the maximum element of $G$ (if such an element exists). See \cite[Proposition 2.10]{HN}.
		\item Note that $\mathfrak{F}_k(M_{\bullet})$ is the smallest quotient in the family $$\{f_j(M_k) \mid j > k\} \cup \{M_k\},$$ that is, the unique object in this set which is a quotient of all the others.	
	\end{enumerate}
\end{obs}}

{\edb Dually, we define a functor  
\[
\mathfrak{T}: \mathcal{N}(\Gamma) \longrightarrow \mathcal{N}^{*}(\Gamma),
\]
such that for any $N_{\bullet} = (N_k)_{k \in G} \in \mathcal{N}(\Gamma)$ we have $\mathfrak{T}(N_{\bullet}) = (\mathfrak{T}_k(N_{\bullet}))_{k \in G}$,  
 where each $\mathfrak{T}_k(N_{\bullet})$ is the torsion part of $N_k$ with respect to the torsion pair
$(\bigcap\limits_{j<k} \mathcal{T}_j, (\bigcap\limits_{j<k} \mathcal{T}_j)^{\bot}) $	
	
In the dual situation,

\[
\mathfrak{T}_k(N_{\bullet}) :=
\begin{cases}
	t_m(N_k), & \text{if $k$ is not the minimum element of $G$,} \\
	N_k, & \text{if $k$ is the minimum element of $G$.}
\end{cases}
\]
where $m$ is any index such that the set $\{t_j(N_k) \mid m \leq j < k\}$ has exactly one element.}

Observe that if $\mathcal{S}'$ is another skeleton of $\mathcal{A}$, and if $\mathfrak{F}'_k$ and $\mathfrak{T}'_k$ denote the functors defined analogously with respect to $\mathcal{S}'$, then $\mathfrak{F}'_k$ and $\mathfrak{T}'_k$ are naturally isomorphic to $\mathfrak{F}_k$ and $\mathfrak{T}_k$, respectively.

We now introduce the functors stratum and substratum, which will serve as the main tools for constructing stratifying systems.

\begin{deff}[Stratum and Substratum]
	\label{stratum}
	Let $\Gamma = \{(\mathcal{T}_k, \mathcal{F}_k) \}_{k \in G}$ be a nested family of torsion pairs in an abelian length category $\mathcal{A}$. 
	\begin{enumerate}
		\item If $M_{\bullet} = (M_k)_{k \in G} \in \mathcal{M}(\Gamma)$, we define the {\bf stratum of $\boldsymbol{M}_{\bullet}$} as 
		\[
		\mathfrak{F}(M_{\bullet}) = (\mathfrak{F}_k(M_{\bullet}))_{k \in G} \in \GR_G(\mathcal{A}).
		\]
		\item If $N_{\bullet} = (N_k)_{k \in G} \in \mathcal{N}(\Gamma)$, we define the {\bf substratum of $\boldsymbol{N}_{\bullet}$} as 
		\[
		\mathfrak{T}(N_{\bullet}) = (\mathfrak{T}_k(N_{\bullet}))_{k \in G} \in \GR_G(\mathcal{A}).
		\]
	\end{enumerate}
\end{deff}

\begin{obs}
	\begin{enumerate} 
		\item Note that both the stratum and the substratum depend on the choice of the nested family $\Gamma = \{(\mathcal{T}_k, \mathcal{F}_k)\}_{k \in G}$.
		\item If $G = \mathbb{I}_n$, then $\mathfrak{T}_k(N_{\bullet}) = t_{k-1}(N_k)$ and $\mathfrak{F}_k(M_{\bullet}) = f_{k+1}(M_k)$, where we set $t_0(N_1) = N_1$ and $f_{n+1}(M_n) = M_n$.
	\end{enumerate}
\end{obs}

We now verify that the constructions introduced above indeed define functors between the intended categories.  
The following lemma ensures that the stratum and the substratum of an object belong to $\mathcal{M}^{*}(\Gamma)$ and $\mathcal{N}^{*}(\Gamma)$, respectively, confirming that  
\[
\mathfrak{F}: \mathcal{M}(\Gamma) \longrightarrow \mathcal{M}^{*}(\Gamma)
\quad \text{and} \quad
\mathfrak{T}: \mathcal{N}(\Gamma) \longrightarrow \mathcal{N}^{*}(\Gamma)
\]
are well-defined functors.

\begin{lema}
	\label{lesma1}
	Let $\Gamma = \{(\mathcal{T}_k, \mathcal{F}_k)\}_{k \in G}$ be a nested family of torsion pairs in an abelian length category $\mathcal{A}$, and let $M_{\bullet} = (M_k)_{k \in G} \in \GR_G(\mathcal{A})$.  
	The following statements hold:
	\begin{enumerate}
		\item If $M_{\bullet} \in \mathcal{M}(\Gamma)$, then $\mathfrak{F}(M_{\bullet}) \in \mathcal{M}^{*}(\Gamma)$. Moreover, $\mathfrak{F}(M_{\bullet}) = M_{\bullet}$ if and only if $M_{\bullet} \in \mathcal{M}^{*}(\Gamma)$.
		\item If $M_{\bullet} \in \mathcal{N}(\Gamma)$, then $\mathfrak{T}(M_{\bullet}) \in \mathcal{N}^{*}(\Gamma)$. Moreover, $\mathfrak{T}(M_{\bullet}) = M_{\bullet}$ if and only if $M_{\bullet} \in \mathcal{N}^{*}(\Gamma)$.
	\end{enumerate}
\end{lema}

\begin{proof}	
We will prove the first item, with the proof of the second being analogous.

1) To show that $\mathfrak{F}(M_{\bullet}) \in \mathcal{M}^{*}(\Gamma)$, we must show for each $k \in G$ that (i) $\mathfrak{F}_k(M_{\bullet}) \neq 0$, (ii) $\mathfrak{F}_k(M_{\bullet}) \in \mathcal{T}_k$ and (iii) $\mathfrak{F}_k(M_{\bullet}) \in \mathcal{F}_j$ if $j > k$ (note that conditions (i) and (iii) imply that $\mathfrak{F}_k(M_{\bullet}) \notin \mathcal{T}_j$ whenever $j > k$).

(i) If $k$ is a maximum element of $G$, then $\mathfrak{F}_k(M_{\bullet}) = M_k \neq 0$. 
Otherwise, by definition, there exists $m > k$ such that $\mathfrak{F}_k(M_{\bullet}) = f_m(M_k)$. 
Since $M_k \notin \mathcal{T}_m$, it follows that $\mathfrak{F}_k(M_{\bullet}) = f_m(M_k) \neq 0$.

(ii) We show that $\mathfrak{F}_k(M_{\bullet}) \in \mathcal{T}_{j}$ whenever $j \leq k$. 
If $k$ is a maximum element of $G$, then $\mathfrak{F}_k(M_{\bullet}) = M_k \in \mathcal{T}_k$, 
and hence $\mathfrak{F}_k(M_{\bullet}) \in \mathcal{T}_j$ for all $j \in G$, 
since the family of torsion pairs is nested.

Suppose otherwise. By definition, there exists $m > k$ such that $\mathfrak{F}_k(M_{\bullet}) = f_m(M_k)$. 
Consider the short exact sequence
\[
0 \longrightarrow t_m(M_k) \longrightarrow M_k \longrightarrow \mathfrak{F}_k(M_{\bullet}) \longrightarrow 0.
\]
Since $\mathfrak{F}_k(M_{\bullet})$ is a quotient of $M_k$ and $M_k \in \mathcal{T}_k \subseteq \mathcal{T}_j$, 
it follows that $\mathfrak{F}_k(M_{\bullet}) \in \mathcal{T}_j$ for all $j \leq k$, 
because the family of torsion pairs is nested.

(iii) We have that $\mathfrak{F}_k(M_{\bullet}) = f_m(M_k) \in \mathcal{F}_j$ for every $j > k$. Indeed, $\mathfrak{F}_k(M_{\bullet}) = f_l(M_k) \in \mathcal{F}_l$ for all $k < l \leq m$, and consequently $\mathfrak{F}_k(M_{\bullet}) \in \mathcal{F}_j$ for $j > k$, since the family of torsion pairs is nested. This shows that $\mathfrak{F}(M_{\bullet}) \in \mathcal{M}^{*}(\Gamma)$.

Finally, since \( \mathfrak{F}(M_{\bullet}) \in \mathcal{M}^{*}(\Gamma) \), we have that \( \mathfrak{F}(M_{\bullet}) = M_{\bullet} \) implies \( M_{\bullet} \in \mathcal{M}^{*}(\Gamma) \). Now, conversely, suppose \( M_{\bullet} \in \mathcal{M}^{*}(\Gamma) \). Then $$ \{f_j(M_k) \mid j > k\} \cup \{M_k\} = \{M_k\}$$ for all $k \in G$, and therefore \( \mathfrak{F}(M_{\bullet}) = M_{\bullet} \).

\end{proof}

\begin{cor}
	\label{cor1}
	Let $\Gamma = \{(\mathcal{T}_k, \mathcal{F}_k)\}_{k \in G}$ be a nested family of torsion pairs in an abelian length category $\mathcal{A}$, and let $M_{\bullet} = (M_k)_{k \in G} \in \GR_G(\mathcal{A})$ with $M_{\bullet} \neq 0_{\bullet}$.
	\begin{enumerate}
		\item If $M_{\bullet}$ is in $\mathcal{M}(\Gamma)$, then for all $k \in G$:
		\begin{enumerate}
			\item[a.] $\mathfrak{F}_k(M_{\bullet}) \neq 0$.
			\item[b.] $\mathfrak{F}_k(M_{\bullet}) \in \mathcal{T}_{j}$ if $ j \leq k$.
			\item[c.] $\mathfrak{F}_k(M_{\bullet}) \in \mathcal{F}_{j}$ if $j > k$.
		\end{enumerate}  
		\item If $M_{\bullet}$ is in $\mathcal{N}(\Gamma)$, then for all $k \in G$:
		\begin{enumerate}
			\item[a'.] $\mathfrak{T}_k(M_{\bullet}) \neq 0$.
			\item[b'.] $\mathfrak{T}_k(M_{\bullet}) \in \mathcal{F}_{j}$ if $ j \geq k$.
			\item[c'.] $\mathfrak{T}_k(M_{\bullet}) \in \mathcal{T}_{j}$ if $j < k$.
		\end{enumerate} 
	\end{enumerate}
\end{cor}
\qed

We conclude this section with an example of how to construct objects in the category $\mathcal{M}(\Gamma)$ and their stratums.

\begin{exem}
	Let $A = kQ$, where $Q$ is the following quiver: 
	\[\begin{tikzcd}
		1 & 2 & 3
		\arrow["\alpha"', from=1-2, to=1-1]
		\arrow["\beta"', from=1-3, to=1-2]
	\end{tikzcd}\] 
	and $\Gamma = \{(\mathcal{T}_k, \mathcal{F}_k)\}_{k \in \mathbb{I}_2}$ be a nested family of torsion pairs, where: 
	\begin{align*}
		\mathcal{T}_1 = \add\{P(1), P(2), S(2)\} \hspace{0.85cm} \mathcal{F}_1 = \add\{S(3)\},\\
		\mathcal{T}_2 = \add \{P(1)\} \hspace{0.9cm} \mathcal{F}_2 =\add \{S(2), S(3), I(2)\}.
	\end{align*}
	
	We describe the non-zero objects in $\mathcal{M}(\Gamma)$. Let $M_{\bullet} = (M_k)_{k \in \mathbb{I}_2} \in \mathcal{M}(\Gamma)$ be a non-zero object, where $M_k \in \mathcal{T}_k$ and $M_k \notin \mathcal{T}_j$ if $j > k$. It is easy to see that $M_1 = P(2)^{\oplus n_1} \bigoplus S(2)^{\oplus n_2}$ and $M_2 = P(1)^{\oplus n_3}$, with $n_1 + n_2 \neq 0$ and $n_3 \neq 0$, since $M_1$ and $M_2$ are non-zero objects. Similarly, it can be shown that the non-zero objects in $\mathcal{M}^{*}(\Gamma)$ are of the form  $M_{\bullet} = (M_k)_{k \in \mathbb{I}_2}$ with $M_1 = S(2)^{\oplus n_1}$ and $M_2 = P(1)^{\oplus n_2}$, with $n_1,n_2 \neq 0$.
	
	We will find the stratum $\mathfrak{F}(M_\bullet)$ of $M_{\bullet} = (M_k)_{k \in \mathbb{I}_2}$. We have that \[\mathfrak{F}_1(M_{\bullet}) = f_2(M_1) = f_2(P(2)^{\oplus n_1} \bigoplus S(2)^{\oplus n_2}) = S(2)^{\oplus (n_1 + n_2)}\] and \[\mathfrak{F}_2(M_{\bullet}) = f_3(M_2) = f_3(P(1)^{\oplus n_3}) = P(1)^{\oplus n_3},\] so 
	\[\mathfrak{F}(M_{\bullet}) = (S(2)^{\oplus (n_1 + n_2)} , P(1)^{\oplus n_3}).\]
	
\end{exem}

\section{Conditions for the existence of a stratifying system}\label{sec4}

\subsection*{\bf Stratifying Systems and Pre-Stratifying Systems}

In this section, we will show that every pre-stratifying system $\Omega = \{\Omega_k\}_{k \in G}$ induces two distinct nested families of torsion pairs $\Gamma$ and $\Gamma'$, both indexed by $G$, such that the characteristic object $M^{\Omega}_{\bullet}$ of $\Omega$ belongs to $\mathcal{M}^{*}(\Gamma)$ and to $\mathcal{N}^{*}(\Gamma')$ (see Definition \ref{deff313}). Conversely, let $\Gamma$ be a nested family of torsion pairs indexed by $G$. If $M_{\bullet}$ is an object in $\mathcal{M}^{\dagger}(\Gamma)$ (or $\mathcal{N}^{\dagger}(\Gamma)$, see Definition \ref{aqui}), then $M_{\bullet}$ induces a pre-stratifying system $\Omega$ such that $M^{\Omega}_{\bullet}$ is a quotient of $M_{\bullet}$ (or a subobject of $M_{\bullet}$, respectively), and consequently, induces at least one stratifying system. We will also show that every stratifying system can be obtained in this way.

In \cite{sistemainfinito}, H. Treffinger demonstrated the existence of stratifying systems of infinite size indexed by $\mathbb{N}$ in $n$-representation infinite algebras. Since there does not always exist an order-preserving bijection between totally ordered sets $(G, \leq)$ and $(H, \preceq)$, even if $G$ and $H$ have the same cardinality, we will consider a more general definition of stratifying systems, maintaining the same philosophy that a stratifying system is a set of objects satisfying certain orthogonality conditions. 

In this direction, we introduce the notion of pre-stratifying systems, which relaxes the indecomposability condition while preserving the orthogonality properties. This broader concept allows the simultaneous construction of several stratifying systems from a single object.

\begin{deff}[Pre-stratifying system and stratifying system]
	\label{definicao}
	Let $(G, \leq)$ be {\edg a totally} ordered set and let $\mathcal{A}$ be an abelian length category.  
	A {\bf pre-stratifying system indexed by $\boldsymbol{G}$} is a collection $\Omega = \{\Omega_k\}_{k \in G}$ of non-zero objects in $\mathcal{A}$ satisfying the following conditions:
	\begin{enumerate}
		\item $\Hom_{\mathcal{A}}(\Omega_k, \Omega_j) = 0$ whenever $k > j$.
		\item $\Ext_{\mathcal{A}}^1(\Omega_k, \Omega_j) = 0$ whenever $k \geq j$.
	\end{enumerate}
	If, in addition, each object $\Omega_k$ is indecomposable, then $\Omega$ is called a {\bf stratifying system}.
\end{deff}

\begin{obs}
	We will denote a pre-stratifying system by $\Delta = \{\Delta_k\}_{k \in G}$ whenever it is, in fact, a stratifying system.
\end{obs}

Let $\Omega = \{\Omega_k\}_{k \in G}$ be a pre-stratifying system in an abelian length category $\mathcal{A}$.
Since $\mathcal{A}$ is a Krull–Schmidt category (as $\mathcal{A}$ is Hom-finite), for each $k \in G$ there exists a decomposition 
\[\Omega_k \cong M_{1} \oplus \cdots \oplus M_{t}\] into indecomposable objects $M_j \in \mathcal{A}$. Such a decomposition is unique up to isomorphism and permutation of the summands. Hence, by choosing an {\edg indecomposable direct summand} $\Delta_k$ of $\Omega_k$ for each $k \in G$, we obtain that the collection $\Delta = \{\Delta_k\}_{k \in G}$ is a stratifying system indexed by $G$. This shows that a pre-stratifying system $\Omega$ yields $\prod_{k \in G} |\Omega_k|$ stratifying systems, where the product should be interpreted as $\infty$ if there are infinitely many indices $k$ such that $|\Omega_k| > 1$.

\subsection*{\bf Filtered modules}
Given a stratifying system $\Delta$ in an abelian length category $\mathcal{A}$, we denote by $\mathcal{F}(\Delta)$ the full subcategory of $\mathcal{A}$ containing the zero object and all objects that are filtered by objects in $\Delta$. That is, a non-zero object $M$ belongs to $\mathcal{F}(\Delta)$ if there exists a finite chain
\[ \eta: 0 = M_0 \subseteq M_1 \subseteq \cdots \subseteq M_n = M \]
of subobjects of $M$ such that $M_i/M_{i-1}$ is isomorphic to an object in $\Delta$ for all $i = 1, 2, \cdots, n$; such a chain is called a {\bf $\boldsymbol{\Delta}$-filtration}.

We say that $\mathcal{F}(\Delta)$ is a {\bf Jordan-H\"older category} if, given $M \in \mathcal{F}(\Delta)$, the multiplicity $[M: \Delta_k]_{\eta}$ does not depend on the filtration $\eta$. Here, $$[M: \Delta_k]_{\eta}$$ denotes the number of quotients isomorphic to $\Delta_k$ in a $\Delta$-filtration $\eta$ (see \cite[Definition 3.1]{jordan} for a more general definition).

K. Erdmann and C. S\'aenz proved that $\mathcal{F}(\Delta)$ is a Jordan-H\"older category if $\Delta$ is a stratifying system of finite size \cite[Lemma 1.4]{definicaosistemaestratificante} in $\modd A$. Using this result, we can show that $\mathcal{F}(\Delta)$ is a Jordan-H\"older category even if $\Delta$ is a stratifying system of infinite size. 

\begin{teo}\label{1234} Let $\Delta = \{\Delta_k\}_{k \in G}$ be a stratifying system in $\modd A$, then $\mathcal{F}(\Delta)$ is a Jordan-H\"older category.    
\end{teo}

\begin{proof} Let $M \in \mathcal{F}(\Delta)$ and consider two $\Delta$-filtrations of $M$:
\[ \eta: 0 = M_0 \subseteq M_1 \subseteq \cdots \subseteq M_n = M \]
and 
\[\eta': 0 = M'_0 \subseteq M'_1 \subseteq \cdots \subseteq M'_m = M. \]
Consider the set $$\mathcal{S} = \{M_i/M_{i-1} \ \mid \ i = 1, 2, \cdots, n \} \cup \{M'_i/M'_{i-1} \ \mid \ i = 1, 2, \cdots, m \}.$$ We have that $\mathcal{S}$ is a finite subset of $\Delta$, thus $\mathcal{S}$ is a stratifying system of finite size with the order inherited from $\Delta$. Furthermore, since $\eta$ and $\eta'$ are $\mathcal{S}$-filtrations, it follows that $n = m$ and $[M: \Delta_k]_{\eta} = [M: \Delta_k]_{\eta'}$ by \cite[Lemma 1.4]{definicaosistemaestratificante}.

\end{proof}

The objects in $\Delta$ are called {\bf relative simple objects of $\boldsymbol{\mathcal{F}(\Delta)}$} because they admit a unique $\Delta$-filtration $\eta: 0 \subseteq \Delta_k$.

\subsection*{\bf Nested Families and Pre-Stratifying Systems}

\begin{deff}[Characteristic object of $\Omega$]\label{deff313}
	Let $\Omega = \{\Omega_k\}_{k \in G}$ be a pre-stratifying system in an abelian length category $\mathcal{A}$ and let $\GR_G(\mathcal{A})$ be the additive category of $G$-graded objects over $\mathcal{A}$. We call the object $ M^{\Omega}_{\bullet} := (\Omega_k)_{k \in G}\in \GR_G(\mathcal{A})$ the {\bf characteristic object of the pre-stratifying system $\boldsymbol \Omega $}.
\end{deff}

Let $M_{\bullet} = (M_k)_{k \in G} \in \GR_G(\mathcal{A})$ and let $\Gamma$ be a nested family of torsion pairs.  
We say that {\bf $M_{\bullet}$ T-induces the nested family $\boldsymbol{\Gamma}$} if  
\[  
\Gamma = \{\,(\mathsf{T}(\{M_j \mid j \ge k\}), (\{M_j \mid j \ge k\})^{\bot})\,\}_{k \in G}.  
\]  
Similarly, we say that {\bf $M_{\bullet}$ F-induces the nested family $\boldsymbol{\Gamma}$} if  
\[  
\Gamma = \{\,({}^{\bot}(\{M_j \mid j \leq k\}), \mathsf{F}(\{M_j \mid j \leq k\}))\,\}_{k \in G}.  
\]  
Note that if $M_{\bullet}$ T-induces $\Gamma$, then $M_{\bullet} \in \mathcal{M}(\Gamma)$; analogously, if $M_{\bullet}$ F-induces $\Gamma$, then $M_{\bullet} \in \mathcal{N}(\Gamma)$. 

We will also say that a pre-stratifying system $\Omega$ {\bf T-induces} (or {\bf F-induces}) the nested family of torsion pairs $\Gamma$ if $M^{\Omega}_{\bullet}$ T-induces (or F-induces, respectively) $\Gamma$.  
For simplicity, we will say that an object {\bf induces} a nested family of torsion pairs whenever it either T-induces or F-induces it, when it is not necessary to distinguish between the two cases.

The following theorem shows how to construct two distinct nested families from an object $ M_{\bullet} = (M_k)_{k \in G}$ such that $ M_k \neq 0$ for all $k\in G$ and $ \Hom_{\mathcal{A}}(M_j, M_i) = 0 $ if $ j > i $. In particular, every pre-stratifying system induces two distinct nested families of torsion pairs through its characteristic object $M^{\Omega}_{\bullet}$.

\begin{teo}\label{galo}
	Let $\mathcal{A}$ be an abelian length category and let $M_{\bullet} = (M_k)_{k \in G}$ be an object of $\GR_G(\mathcal{A})$ with $M_k\neq 0$ for all $k\in G$ such that $\Hom_{\mathcal{A}}(M_j, M_i) = 0 $ if $j > i$. Then the following statements hold true.

	\begin{enumerate}
		\item $M_{\bullet}$ T-induces the nested family $$\Gamma = \{\,(\mathsf{T}(\{M_j \mid j \ge k\}), (\{M_j \mid j \ge k\})^{\bot})\,\}_{k \in G}.$$
		\item $M_{\bullet}$ F-induces the nested family $$\Gamma' = \{\,({}^{\bot}(\{M_j \mid j \leq k\}), \mathsf{F}(\{M_j \mid j \leq k\}))\,\}_{k \in G}.$$
		\item $\Gamma$ and $\Gamma'$ are distinct nested families of torsion pairs.
		\item The object $M_{\bullet}$ lies in the intersection $\mathcal{M}^{*}(\Gamma) \bigcap \mathcal{N}^{*}(\Gamma')$.
	\end{enumerate}
\end{teo}

\begin{proof}
1) Since $\{M_j \mid j \ge i\} \subsetneq \{M_j \mid j \ge k\}$ if $i > k$, then $$\mathsf{T}(\{M_j \mid j \ge i\}) \subseteq \mathsf{T}(\{M_j \mid j \ge k\}).$$

\noindent The class $\mathsf{T}(\{M_j \mid j \ge k\})$ is defined as the smallest torsion class containing the objects in $\{M_j \mid j \ge k\}$, so $M_k \in \mathsf{T}(\{M_j \mid j \ge k\})$. On the other hand, $M_k \in (\{M_j \mid j \ge i\})^{\bot}$ whenever $i > k$. Hence $M_k \notin \mathsf{T}(\{M_j \mid j \ge i\})$ for all $i > k$. Consequently, $\Gamma$ is a nested family of torsion pairs. 

2) This follows by a similar argument to 1).

3) It suffices to note that $$M_k \in \mathsf{T}(\{M_j \mid j \ge k\})$$ and $$M_k \notin {}^{\bot}(\{M_j \mid j \leq k\})$$ for all $k \in G$, hence $\mathsf{T}(\{M_j \mid j \ge k\}) \neq {}^{\bot}(\{M_j \mid j \leq k\})$, and consequently $\Gamma \neq \Gamma'$.

4) We have that $M_{\bullet}\in \mathcal{M}(\Gamma)$ since, for all $k \in G$, $$M_k \neq 0, M_k \in \mathsf{T}(\{M_j \mid j \ge k\}) \quad \hbox{ and } \quad M_k \notin \mathsf{T}(\{M_j \mid j \ge i\})$$ for all $i > k$. On the other hand, since $M_k \in (\{M_j \mid j \ge i\})^{\bot}$ whenever $i > k$, it follows that $M_{\bullet}\in \mathcal{M}^{*}(\Gamma)$. Similarly, $M_{\bullet}\in \mathcal{N}^{*}(\Gamma')$, which completes the proof.
\end{proof}

\begin{cor} \label{nestedfamilyinducedbystrat}
	Let $ \Omega = \{\Omega_k\}_{k \in G} $ be a pre-stratifying system in an abelian length category $\mathcal{A}$. Then $ \Omega $ induces two distinct nested families of torsion pairs $ \Gamma $ and $ \Gamma' $, both indexed by $ G $, such that $ M^{\Omega}_{\bullet} $ is in $ \mathcal{M}^{*}(\Gamma) $ and in $ \mathcal{N}^{*}(\Gamma') $.
\end{cor}
\qed

The following example illustrates how to obtain a non-enumerable nested family that cannot be induced from a stratifying system.

\begin{exem}\label{exem2}
	Let $ A = \mathbb{C}Q$ be the Kronecker algebra, and let $ (\mathbb{C}, \leq)$ be the field of complex numbers with the lexicographic order, that is, $ x + yi > x' + y'i $ if $ x > x' $ or $ x = x' $ and $ y > y' $. For every $ \lambda \in \mathbb{C} $, consider the module $ M_{\lambda} $ whose representation is given by
	\[
	\begin{tikzcd}
		\mathbb{C} & \mathbb{C}
		\arrow["\lambda"', shift right, from=1-2, to=1-1]
		\arrow["1", shift left, from=1-2, to=1-1]
	\end{tikzcd}
	\]
	We have that the object $M_{\bullet} = ( M_{\lambda})_{\lambda \in \mathbb{C}} \in \GR_{\mathbb{C}}(\modd A)$ satisfies the conditions of the Theorem \ref{galo} (since $\Hom_A(M_{\lambda},M_{\rho}) = 0$ if $\lambda \neq \rho$), thus inducing two uncountable nested families of torsion pairs. However, {\edg there are no} stratifying systems of size greater than 2 in the Kronecker algebra \cite[Theorem 2.3.18]{tesecadavid}.
\end{exem}

Let $\Gamma$ be a nested family of torsion pairs in an abelian length category $\mathcal{A}$, and let $M_{\bullet} = (M_k)_{k \in G} \in \mathcal{M}(\Gamma)$. Corollary \ref{cor1} ensures that $\mathfrak{F}(M_{\bullet})$ satisfies $\Hom_{\mathcal{A}}(\mathfrak{F}_k(M_{\bullet}), \mathfrak{F}_j(M_{\bullet})) = 0$ whenever $k > j$. If, in addition, one can ensure that $\Ext_{\mathcal{A}}^{1}(\mathfrak{F}_k(M_{\bullet}), \mathfrak{F}_j(M_{\bullet})) = 0$ for $k \ge j$, then $\Omega = \{\mathfrak{F}_k(M_{\bullet})\}_{k \in G}$ is a pre-stratifying system. In particular, by choosing an {\edg indecomposable direct summand} $\Delta_k$ of $\mathfrak{F}_k(M_{\bullet})$ for each $k \in G$, we obtain a stratifying system $\Delta = \{\Delta_k\}_{k \in G}$. This observation motivates the following definition.

\begin{deff}
	\label{aqui}
	Let $\Gamma = \{(\mathcal{T}_k, \mathcal{F}_k) \}_{k \in G}$ be a nested family of torsion pairs in an abelian length category $\mathcal{A}$.  
	We define $\mathcal{M}^{\dagger}(\Gamma)$ as the full subcategory of $\mathcal{M}(\Gamma)$ consisting of all non-zero objects $M_{\bullet} = (M_k)_{k \in G}$ satisfying $\Ext_{\mathcal{A}}^{1}(M_k, \mathfrak{F}_j(M_{\bullet})) = 0$ for all $k \ge j$.
	
	Similarly, we define $\mathcal{N}^{\dagger}(\Gamma)$ as the full subcategory of $\mathcal{N}(\Gamma)$ consisting of all non-zero objects $N_{\bullet} = (N_k)_{k \in G}$ satisfying $\Ext_{\mathcal{A}}^{1}(\mathfrak{T}_k(N_{\bullet}), N_j) = 0$ for all $k \ge j$.  
\end{deff}

Let $\Gamma = \{(\mathcal{T}_k, \mathcal{F}_k)\}_{k \in G}$ be a nested family of torsion pairs in an abelian length category $\mathcal{A}$, $\Omega = \{\Omega_k\}_{k \in G}$ be a pre-stratifying system in $\mathcal{A}$, and let $M_{\bullet}$ be an object in $\mathcal{M}(\Gamma)$. We say that {\bf $\boldsymbol{M_{\bullet}}$ induces the pre-stratifying system $\boldsymbol{\Omega = \{\Omega_k\}_{k \in G}}$ as a quotient} if $\Omega_k = \mathfrak{F}_k(M_{\bullet})$ for all $k \in G$. Similarly, for $N_{\bullet}$ in $\mathcal{N}(\Gamma)$, we say that {\bf $\boldsymbol{N_{\bullet}}$ induces the pre-stratifying system $\boldsymbol{\Omega = \{\Omega_k\}_{k \in G}}$ as a subobject} if $\Omega_k = \mathfrak{T}_k(N_{\bullet})$ for all $k \in G$. 

The following theorem shows that objects in $\mathcal{M}^{\dagger}(\Gamma)$ (respectively, in $\mathcal{N}^{\dagger}(\Gamma)$) induce a pre-stratifying system as a quotient (respectively, as a subobject).

\begin{teo}
	\label{teoresma1}
	Let $\Gamma = \{(\mathcal{T}_k, \mathcal{F}_k) \}_{k \in G}$ be a nested family of torsion pairs in an abelian length category $\mathcal{A}$ and let $M_{\bullet} = (M_k)_{k \in G}$ be an object in $\GR_G(\mathcal{A})$. 
	
	\begin{enumerate}
		\item If $M_{\bullet} \in \mathcal{M}^{\dagger}(\Gamma)$, then $M_{\bullet}$ induces the pre-stratifying system $\Omega = \{\mathfrak{F}_k(M_{\bullet})\}_{k \in G}$ as a quotient.
		
		\item If $M_{\bullet} \in \mathcal{N}^{\dagger}(\Gamma) $, then $M_{\bullet}$ induces the pre-stratifying system $\Omega = \{\mathfrak{T}_k(M_{\bullet})\}_{k \in G}$ as a subobject.
	\end{enumerate}
\end{teo}

\begin{proof} 1) Assume that \( M_{\bullet} \in \mathcal{M}^{\dagger}(\Gamma) \). By Corollary~\ref{cor1}, it follows that \( \mathfrak{F}_k(M_{\bullet}) \neq 0 \) for every \( k \in G \).
	
	Suppose first that $i > j$. Then $\mathfrak{F}_i(M_{\bullet}) \in \mathcal{T}_i$ 
	and $\mathfrak{F}_j(M_{\bullet}) \in \mathcal{F}_i$ by Corollary~\ref{cor1}, 
	hence $\Hom_{\mathcal{A}}(\mathfrak{F}_i(M_{\bullet}), \mathfrak{F}_j(M_{\bullet})) = 0$.

Now assume that $i \geq j$ and that $i$ is not the maximum element of $G$ (if such an element exists). In this case, there exists $m > i \geq j$ such that $\mathfrak{F}_i(M_{\bullet}) = f_m(M_i)$, and we have a short exact sequence
\begin{equation*}
	0 \longrightarrow t_m(M_i) \longrightarrow M_i \longrightarrow \mathfrak{F}_i(M_{\bullet}) \longrightarrow 0.
\end{equation*}
Applying the functor $\Hom_{\mathcal{A}}(-,\mathfrak{F}_j(M_{\bullet}))$ yields
\begin{equation*}
	\Hom_{\mathcal{A}}(t_m(M_i), \mathfrak{F}_j(M_{\bullet})) \longrightarrow \Ext_{\mathcal{A}}^1(\mathfrak{F}_i(M_{\bullet}), \mathfrak{F}_j(M_{\bullet})) \longrightarrow \Ext_{\mathcal{A}}^1(M_i, \mathfrak{F}_j(M_{\bullet})).
\end{equation*}
Since $M_{\bullet} \in \mathcal{M}^{\dagger}(\Gamma)$, we have $\Ext_{\mathcal{A}}^1(M_i, \mathfrak{F}_j(M_{\bullet})) = 0$. Moreover, $t_m(M_i) \in \mathcal{T}_m$, and by Corollary \ref{cor1}, we have that $\mathfrak{F}_j(M_{\bullet}) \in \mathcal{F}_m$, hence $$\Hom_{\mathcal{A}}(t_m(M_i), \mathfrak{F}_j(M_{\bullet})) = 0.$$ Therefore, $\Ext_{\mathcal{A}}^1(\mathfrak{F}_i(M_{\bullet}), \mathfrak{F}_j(M_{\bullet})) = 0$.

If $i$ is the maximum element of $G$, the result follows by the same argument, starting from the short exact sequence
\begin{equation*}
	0 \longrightarrow 0 \longrightarrow M_i \longrightarrow \mathfrak{F}_i(M_{\bullet}) \longrightarrow 0.
\end{equation*}
2) The proof is analogous to that of 1).
\end{proof}

Example \ref{exem2} shows that $\mathcal{M}^{\dagger}(\Gamma)$ may be empty; otherwise, it would be possible to obtain a stratifying system of infinite size in the Kronecker algebra. 

The following theorem provides the converse direction of the previous result, completing the characterization of pre-stratifying systems in abelian length categories. It shows that every such system arises naturally from a suitable nested family of torsion pairs. More precisely, for any pre-stratifying system $\Omega = \{\Omega_k\}_{k \in G}$, there exists a nested family $\Gamma = \{(\mathcal{T}_k, \mathcal{F}_k)\}_{k \in G}$ and an object $M_{\bullet} \in \mathcal{M}^{\dagger}(\Gamma)$ (or $N_{\bullet} \in \mathcal{N}^{\dagger}(\Gamma)$) such that $\Omega$ is induced as a quotient of $M_{\bullet}$ (or, respectively, as a subobject of $N_{\bullet}$).
In particular, since every stratifying system is a pre-stratifying one, it follows that every stratifying system can also be obtained in this way.

\begin{teo}\label{cubo}
	Let $\Omega = \{\Omega_k\}_{k \in G}$ be a pre-stratifying system in an abelian length category $\mathcal{A}$.
	\begin{enumerate}
		\item There exists at least one nested family of torsion pairs $$\Gamma = \{(\mathcal{T}_k, \mathcal{F}_k)\}_{k \in G}$$ and at least one object $M_{\bullet}$ in $\mathcal{M}^{\dagger}(\Gamma)$ such that $\Omega$ is induced as a quotient of $M_{\bullet}$.
		\item There exists at least one nested family of torsion pairs $$\Gamma' = \{(\mathcal{T}'_k, \mathcal{F}'_k)\}_{k \in G}$$ and at least one object $N_{\bullet}$ in $\mathcal{N}^{\dagger}(\Gamma')$ such that $\Omega$ is induced as a subobject of $N_{\bullet}$.
	\end{enumerate}
\end{teo}

\begin{proof} 1) Suppose that $\Omega = \{\Omega_k\}_{k \in G}$ is a pre-stratifying system, and let $M^{\Omega}_{\bullet} = (\Omega_k)_{k \in G}$ be its characteristic object. Theorem \ref{galo} ensures that 
	\[
	\Gamma = \{\,(\mathsf{T}(\{\Omega_j \mid j \ge k\}), (\{\Omega_j \mid j \ge k\})^{\bot})\,\}_{k \in G}
	\]
	is a nested family of torsion pairs and that $M^{\Omega}_{\bullet} \in \mathcal{M}^{*}(\Gamma)$.
	
	We now show that, in fact, $M^{\Omega}_{\bullet} \in \mathcal{M}^{\dagger}(\Gamma)$. Combining Theorem \ref{galo} with Lemma \ref{lesma1}, we obtain $\mathfrak{F}_k(M^{\Omega}_{\bullet}) = \Omega_k$. Hence,
	
	\[
	\Ext_{\mathcal{A}}^1(\Omega_j, \mathfrak{F}_k(M^{\Omega}_{\bullet})) 
	= \Ext_{\mathcal{A}}^1(\Omega_j, \Omega_k) = 0
	\]
	whenever $j \geq k$, since $\Omega$ is a pre-stratifying system. Therefore, by Theorem \ref{teoresma1}, the family $\Gamma$ together with the object $M^{\Omega}_{\bullet}$ induces the pre-stratifying system $\Omega$. Thus, $M^{\Omega}_{\bullet}$ is one possible choice for the object $M_{\bullet}$ required in (1).
	
	2) This follows by a similar argument to 1).

\end{proof}

\section{Stratifying system via $\tau$-rigid modules}\label{sec5}

In \cite[Theorem 3.4]{stratifying}, H. Treffinger and O. Mendoza showed that every non-zero basic $\tau$-rigid module induces a stratifying system. In this section, we will demonstrate that this result can be obtained as a corollary of Theorem \ref{teoresma1}.

We begin with a lemma providing a sufficient condition for an object $M_{\bullet} \in \mathcal{M}(\Gamma)$ to be in $\mathcal{M}^{\dagger}(\Gamma)$.

\begin{lema}
	\label{extproj}
Let $\Gamma = \{(\mathcal{T}_k, \mathcal{F}_k)\}_{k \in \mathbb{I}_n}$ be a nested family of torsion pairs in $\mathcal{A}$, and let $M_{\bullet} \in \GR_{\mathbb{I}_n}(\mathcal{A})$. Define $M = \bigoplus\limits_{k = 1}^{n} M_k$.
	
	\begin{enumerate}
		\item If $M_{\bullet} \in \mathcal{M}(\Gamma)$ and $M$ is Ext-projective in $\mathcal{T}_1$, then $M_{\bullet} \in \mathcal{M}^{\dagger}(\Gamma)$.
		\item If $M_{\bullet} \in \mathcal{N}(\Gamma)$ and $M$ is Ext-injective in $\mathcal{F}_n$, then $M_{\bullet} \in \mathcal{N}^{\dagger}(\Gamma)$.
	\end{enumerate}
\end{lema}

\begin{proof} 1) Since $M_k$ is a direct summand of $M$, it is Ext-projective in $\mathcal{T}_1$. On the other hand, Corollary \ref{cor1} ensures that $\mathfrak{F}_j(M_{\bullet}) \in \mathcal{T}_1$ for all $j \in \mathbb{I}_n$. Thus, $\Ext_{\mathcal{A}}^1(M_k, \mathfrak{F}_j(M_{\bullet})) = 0$ for all $j \in \mathbb{I}_n$. In particular, $\Ext_{\mathcal{A}}^1(M_k, \mathfrak{F}_j(M_{\bullet})) = 0$ if $k \geq j$ and $M_{\bullet} \in \mathcal{M}^{\dagger}(\Gamma)$. 

2) This follows by a similar argument to 1). 

\end{proof}

The following lemma shows how to induce a nested family of torsion pairs from a certain class of objects.

\begin{lema}\label{TFad}
Let $\mathcal{A}$ be an abelian length category and let $M_{\bullet} = (M_k)_{k\in G} \in \GR_G(\mathcal{A})$ be such that $M_k$ is a non-zero object for each $k \in G$.
	\begin{enumerate}
		\item If for every $k\in G$, we have $M_k \notin \mathsf{T}(\{M_j \mid j > k\})$, then $M_{\bullet}$ T-induces the nested family $\Gamma = \{\,(\mathsf{T}(\{M_j \mid j \ge k\}), (\{M_j \mid j \ge k\})^{\bot})\,\}_{k \in G}$.
		\item If for every $k \in G$, we have $M_k \notin \mathsf{F}(\{M_j \mid j <k\})$, then $M_{\bullet}$ F-induces the nested family $\Gamma = \{\,({}^{\bot}(\{M_j \mid j \leq k\}), \mathsf{F}(\{M_j \mid j \leq k\}))\,\}_{k \in G}$.
	\end{enumerate}
\end{lema}

\begin{proof} 1) Clearly, $\mathsf{T}(\{M_j \mid j > k\}) \subseteq \mathsf{T}(\{M_j \mid j \ge k\})$ and $$M_k \in \mathsf{T}(\{M_j \mid j \ge k\}).$$ Since $M_k \notin \mathsf{T}(\{M_j \mid j > k\})$, it follows that $\Gamma$ is a nested family of torsion pairs.

2) The second case is analogous. 

\end{proof}

Unlike in Theorem \ref{galo}, it is not necessary that $M_{\bullet}$ is in $\mathcal{M}^{*}(\Gamma)$ (or in $\mathcal{N}^{*}(\Gamma)$) if $M_{\bullet}$ induces $\Gamma$ as in the lemma above.

\begin{exem}\label{exem3}  
	Let $P(1), P(2), \ldots, P(n)$ be an enumeration of the indecomposable projective modules in $\modd A$.  
	Given a permutation $\sigma: \mathbb{I}_n \to \mathbb{I}_n$ and a subset $\{k_2, k_3, \ldots, k_t\} \subseteq \mathbb{I}_{n-1}$ with $t-1$ elements such that $$k_2 < k_3 < \cdots < k_t,$$ define a decomposition $A_A = \bigoplus\limits_{j=1}^{t} M_j$, where $M_j = \bigoplus\limits_{i = k_j + 1}^{k_{j+1}} P(\sigma(i))$ for $1 \le j < t$, and  
	$M_t = \bigoplus\limits_{i = k_t + 1}^{n} P(\sigma(i))$ (here we fix $k_1 = 0$ and $\mathbb{I}_0 = \{\}$).  
	
	It is easy to see that $M_j \notin \Fac\!\left(\bigoplus\limits_{k > j}^{t} M_k\right)$, since every epimorphism onto $M_j$ splits. Therefore, $M_{\bullet} = (M_k)_{k \in \mathbb{I}_t}$ induces a nested family indexed by $\mathbb{I}_t$.  
	This shows that any decomposition of $A_A$ is torsion-free admissible, and that $A_A$ induces $2^{n-1}n!$ distinct nested families.  
	A similar argument shows that $D({}_A A)$ induces $2^{n-1}n!$ distinct nested families as well.  
\end{exem}

Let $M$ be a $\tau$-rigid module in $\modd A$. Proposition \ref{taurigidoalmostsplit} ensures that $\mathsf{T}(M) = \Fac(M)$ since it always holds that $\Fac(M) \subseteq \mathsf{T}(M)$. Using mutation techniques from $\tau$-tilting theory, H. Treffinger and O. Mendoza proved that every basic non-zero $\tau$-rigid module admits at least one torsion-free admissible indecomposable decomposition \cite[Proposition 3.2]{stratifying}. Assuming the existence of such a decomposition, we will show that \cite[Theorem 3.4]{stratifying} can be seen as a corollary of Theorem \ref{teoresma1}.

{\edb \begin{cor}\label{teoremadoscara}
	Let $M$ be a non-zero, basic, and $\tau$-rigid module with a torsion-free admissible decomposition $M = \bigoplus\limits_{j=1}^{t}M_j$.	
	
	\begin{enumerate}
		\item Let $r \leq t$ and consider a subset $\{k_2, k_3, \ldots, k_r\} \subseteq \mathbb{I}_{t-1}$ such that $$k_2 < \cdots < k_r.$$ Define a decomposition $M = \bigoplus\limits_{j=1}^{r} L_j$, where $L_j = \bigoplus\limits_{i = k_j + 1}^{k_{j+1}} M_i$ for $1 \le j < r$, and $L_r = \bigoplus\limits_{i = k_r + 1}^{t} M_i$ (here we set $k_1 = 0$ and $\mathbb{I}_0 = \{\}$). 
		
		Then the object $L_{\bullet} = (L_k)_{k \in \mathbb{I}_{r}}$ T-induces the nested family $$\Gamma = \{(\Fac(\bigoplus\limits_{j\geq k}L_j), (\bigoplus\limits_{j\geq k}L_j)^{\bot})\}_{k \in \mathbb{I}_{r}},$$ $L_{\bullet} \in \mathcal{M}^{\dagger}(\Gamma)$ and $\Omega = \{\Omega_k := f_{k+1}(L_k)\}_{k \in \mathbb{I}_r}$ is a pre-stratifying system of size $r$, where $f_k$ is the torsion-free functor associated to the pair $(\Fac(\bigoplus\limits_{j\geq k}L_j), (\bigoplus\limits_{j\geq k}L_j)^{\bot})$.
		\item \cite[Theorem 3.4]{stratifying} If $f_k$ denotes the torsion-free functor associated to the pair $(\Fac(\bigoplus\limits_{j\geq k}M_j), (\bigoplus\limits_{j\geq k}M_j)^{\bot})$, then $\Delta = \{\Delta_k := f_{k+1}(M_k)\}_{k \in \mathbb{I}_t}$ is a stratifying system of size $t$.
	\end{enumerate}	
\end{cor}

\begin{proof} 1) For $j < r$, we have that $L_j \notin \mathsf{T}(\{L_i \mid i > j\}) = \Fac(\bigoplus\limits_{i > j}L_i)$; otherwise, $M_{k_{j+1}} \in \Fac(\bigoplus\limits_{i > j}L_i) = \Fac(\bigoplus\limits_{i > k_{j+1}}M_i)$, contradicting the fact that the decomposition $M= \bigoplus\limits_{j=1}^{t}M_j$ is torsion-free admissible. Hence, Lemma \ref{TFad} ensures that $L_{\bullet} = (L_k)_{k \in \mathbb{I}_r}$ T-induces $\Gamma = \{(\Fac(\bigoplus\limits_{j\geq k}L_j), (\bigoplus\limits_{j\geq k}L_j)^{\bot})\}_{k \in \mathbb{I}_{r}}$, and consequently $L_{\bullet} \in \mathcal{M}(\Gamma)$. 
	
On the other hand, Proposition \ref{rigidoextproj} and Lemma \ref{extproj} guarantee that $L_{\bullet}\in\mathcal{M}^{\dagger}(\Gamma)$. Therefore, Theorem \ref{teoresma1} implies that $$\Omega = \{\mathfrak{F}_k(L_{\bullet}) = f_{k+1}(L_k)\}_{k \in \mathbb{I}_r}$$ is a pre-stratifying system.

2) Taking $r=t$ in the construction above shows that $$\Delta = \{\Delta_k = f_{k+1}(M_k)\}_{k \in \mathbb{I}_t}$$ is a pre-stratifying system. Moreover, Lemma \ref{ind} ensures that $f_{i+1}(M_i)$ is indecomposable since $M_i \bigoplus \left(\bigoplus\limits_{k > i}M_k \right)$ is $\tau$-rigid and $M_i \notin \Fac(\bigoplus\limits_{k > i}M_k)$; hence, $\Delta$ is a stratifying system.

\end{proof}}

It is possible to prove the dual of Corollary \ref{teoremadoscara} and show that every basic $\tau^{-}$-rigid module (including injective modules) induces at least one stratifying system as a subobject.

\section{Order Relations on Nested Families}\label{sec6}

{\edg Let $M = \bigoplus\limits_{j=1}^{t} M_j$ be a torsion-free admissible decomposition into {\edg indecomposable direct summand}s of the basic $\tau$-rigid module $M$. Corollary \ref{teoremadoscara} shows that $M_{\bullet} = (M_k)_{k \in \mathbb{I}_t} \in \mathcal{M}^{\dagger}(\Gamma)$, where $\Gamma = \{(\Fac(\bigoplus\limits_{j \geq k} M_j), (\bigoplus\limits_{j \geq k} M_j)^{\bot})\}_{k \in \mathbb{I}_t}$. However, it is possible that $M_{\bullet} \in \mathcal{M}^{\dagger}(\Gamma')$ for some $\Gamma' \neq \Gamma$. In this section, we study how the induced pre-stratifying systems relate to each other under different nested families. We begin with a definition.}

\begin{deff}
	Let $\Gamma = \{(\mathcal{T}_k, \mathcal{F}_k)\}_{k \in G}$ and $\Gamma' = \{(\mathcal{T}'_k, \mathcal{F}'_k)\}_{k \in G}$ be nested families of torsion pairs. We say that {\bf $\boldsymbol{\Gamma}$ is coarser than $\boldsymbol{\Gamma'}$}, or equivalently, that {\bf $\boldsymbol{\Gamma'}$ is finer than $\boldsymbol{\Gamma}$}, if $\mathcal{T}'_k \subseteq \mathcal{T}_k$ for all $k \in G$.
\end{deff}

We write $\Gamma' \leq \Gamma$ to indicate that $\Gamma$ is coarser than $\Gamma'$ (or equivalently, that $\Gamma'$ is finer than $\Gamma$). Note that, by definition, every nested family $\Gamma$ is coarser than itself.

The following proposition describes the relationship between the nested families induced by $M_{\bullet}$ and the nested families $\Gamma'$ such that $M_{\bullet} \in \mathcal{M}(\Gamma')$.

\begin{prop}
	Let $\mathcal{A}$ be an abelian length category, $M_{\bullet} = (M_k)_{k \in G}$ an object of $\GR_G(\mathcal{A})$ with each $M_k$ being a non-zero object, and let $\Gamma' = \{(\mathcal{T}_k, \mathcal{F}_k)\}_{k \in G}$ be a nested family of torsion pairs.
	
	\begin{enumerate}
		\item If $M_{\bullet}$ T-induces the nested family $$\Gamma = \{\,(\mathsf{T}(\{M_j \mid j \ge k\}), (\{M_j \mid j \ge k\})^{\bot})\,\}_{k \in G},$$ then $\Gamma$ is the {\bf finest} nested family satisfying $M_{\bullet} \in \mathcal{M}(\Gamma)$, that is, if $M_{\bullet} \in \mathcal{M}(\Gamma')$, then $\Gamma \leq \Gamma'$.
		
		\item If $M_{\bullet}$ F-induces the nested family $$\Gamma = \{\,({}^{\bot}(\{M_j \mid j \leq k\}), \mathsf{F}(\{M_j \mid j \leq k\}))\,\}_{k \in G},$$ then $\Gamma$ is the {\bf coarsest} nested family satisfying $M_{\bullet} \in \mathcal{N}(\Gamma)$, that is, if $M_{\bullet} \in \mathcal{N}(\Gamma')$, then $\Gamma' \leq \Gamma$.
	\end{enumerate}
\end{prop}

\begin{proof} 1) Suppose that $M_{\bullet} \in \mathcal{M}(\Gamma')$. Then for all $k \in G$, we have that $\{M_j \mid j \ge k\} \subseteq \mathcal{T}'_k$. Since $\mathsf{T}(\{M_j \mid j \ge k\})$ is the smallest torsion class containing $\{M_j \mid j \ge k\}$, we have $\mathsf{T}(\{M_j \mid j \ge k\}) \subseteq \mathcal{T}'_k$, hence $\Gamma \leq \Gamma'$.

2) The second case is analogous. 
\end{proof}

Let $M$ be a non-zero, basic, and $\tau$-rigid module with a torsion-free admissible decomposition $M = \bigoplus\limits_{j=1}^{t} M_j$ and consider $M_{\bullet} = (M_k)_{k \in \mathbb{I}_t}$. Corollary \ref{teoremadoscara} shows that $\mathfrak{F}_k(M_{\bullet})$ is a basic module for all $k \in \mathbb{I}_t$, where $\Gamma = \{(\Fac(\bigoplus\limits_{j\geq k}M_j), (\bigoplus\limits_{j\geq k}M_j)^{\bot})\}_{k \in \mathbb{I}_{t}}$. The following proposition shows the relationship between the stratums of $M_{\bullet}$ when changing the family $\Gamma$ to a coarser one.

\begin{prop}\label{final} Let $\mathcal{A}$ be an abelian length category, and let $$\Gamma = \{(\mathcal{T}_k, \mathcal{F}_k)\}_{k \in G} \quad \hbox{ and } \quad \Gamma' = \{(\mathcal{T}'_k, \mathcal{F}'_k)\}_{k \in G}$$ be nested families of torsion pairs in $\mathcal{A}$ with $\Gamma \leq \Gamma'$. Consider an object $M_{\bullet} = (M_k)_{k \in G}$ in $\GR_G(\mathcal{A})$. For all $k \in G$:
	\begin{enumerate}
		\item If $M_{\bullet} \in\mathcal{M}(\Gamma)$ and $M_{\bullet} \in\mathcal{M}(\Gamma')$, then $\mathfrak{F}'_k(M_{\bullet})$ is a quotient of $\mathfrak{F}_k(M_{\bullet})$, where $\mathfrak{F}_k(M_{\bullet})$ and $\mathfrak{F}'_k(M_{\bullet})$ are the $k$-th components of the stratum of $M_{\bullet}$ with respect to $\Gamma$ and $\Gamma'$, respectively.
		\item If $M_{\bullet}\in\mathcal{N}(\Gamma)$ and $M_{\bullet}\in\mathcal{N}(\Gamma')$, then $\mathfrak{T}_k(M_{\bullet})$ is a subobject of $\mathfrak{T}_k'(M_{\bullet})$, where $\mathfrak{T}_k(M_{\bullet})$ and $\mathfrak{T}'_k(M_{\bullet})$ are the $k$-th components of the substratum of $M_{\bullet}$ with respect to $\Gamma$ and $\Gamma'$, respectively.
	\end{enumerate}
\end{prop}

\begin{proof} 
	1) We have that $\mathfrak{F}_k(M_{\bullet})$ and $\mathfrak{F}'_k(M_{\bullet})$ are the smallest quotients of the families $\{f_j(M_k) \mid j > k \} \cup \{M_k\}$ and $\{f'_j(M_k) \mid j > k \} \cup \{M_k\}$, respectively. Since $\mathcal{T}_j \subseteq \mathcal{T}'_j$, the fact that $t_j(M_k)$ is a subobject of $t'_j(M_k)$ implies that $f'_j(M_k)$ is a quotient of $f_j(M_k)$, and then $\mathfrak{F}'_k(M_{\bullet})$ is a quotient of $\mathfrak{F}_k(M_{\bullet})$. 
	
	2) The second case is analogous. 
\end{proof}

The next two examples explore the consequences of Proposition \ref{final}.

\begin{exem}\label{exem101}
	\setlength{\arraycolsep}{-0.15pt}
	\def\arraystretch{0.15} 
	Let $A = kQ/I$, where $Q$ is the quiver
	\[\begin{tikzcd}
		2 & 4 & 6 \\
		1 & 3 & 5
		\arrow["\mu"', from=1-1, to=2-1]
		\arrow["\lambda"', from=1-2, to=1-1]
		\arrow["\beta"', from=1-2, to=2-2]
		\arrow["\alpha"', from=1-3, to=1-2]
		\arrow["\gamma", from=1-3, to=2-3]
		\arrow["\nu", from=2-2, to=2-1]
		\arrow["\delta", from=2-3, to=2-2]
	\end{tikzcd}\]
	and $I$ is the ideal generated by $\alpha\beta - \gamma\delta$ and $\lambda \mu - \beta \nu$. 
	
	Consider the $\tau$-rigid $M = M_1 \bigoplus M_2$, where $M_1 = \begin{matrix}
		2& & 3\\ &1&
	\end{matrix}$ and $M_2 = 6$ (see \cite[Example IV.1.13]{coelho} for the Auslander-Reiten quiver of the algebra $A$ and \cite[Example I.2.18]{coelho} for a brief explanation of the notation). Such a decomposition of $M$ is torsion-free admissible. 
	
	Let $$\Gamma = \{(\Fac(M_1\bigoplus M_2), (M_1 \bigoplus M_2)^{\bot}), (\Fac(M_2), M_2^{\bot})\}.$$
	 Corollary \ref{teoremadoscara} ensures that $(\{\mathfrak{F}_1(M_{\bullet}) = M_1, \mathfrak{F}_2(M_{\bullet}) = M_2\}, \leq)$ is a stratifying system, where $M_{\bullet} = (M_j)_{j \in \mathbb{I}_2}$. 
	
	On the other hand, for the nested family $$\Gamma' = \{(\modd A, \{0\}), (\Fac(1 \bigoplus 6), (1 \bigoplus 6)^{\bot})\},$$ we have $\Gamma \leq \Gamma'$ and $M_{\bullet} \in \mathcal{M}^{\dagger}(\Gamma')$. Moreover, $\mathfrak{F}'_1(M_{\bullet}) = 2 \bigoplus 3$ and $\mathfrak{F}'_2(M_{\bullet}) = 6$. Theorem \ref{teoresma1} ensures that $(\{M'_1 = 2 \bigoplus 3, M'_2 = 6 \}, \leq)$ is a pre-stratifying system, from which we obtain two stratifying systems induced by $M_{\bullet}$: $$(\{\Delta_1 = 2, \Delta_2 = 6\},\leq) \hbox{ and } (\{\Delta'_1=3, \Delta'_2=6\}, \leq).$$
\end{exem}

The above example suggests that an object $M_{\bullet} \in \mathcal{M}^{\dagger}(\Gamma)$ induces a unique stratifying system if $\Gamma$ is a ``sufficiently fine'' nested family. The next example shows that this is not the case.

\begin{exem}
	\setlength{\arraycolsep}{-0.15pt}
	\def\arraystretch{0.15} 
	Let $A$ be the same algebra as in Example \ref{exem101} and consider also the module $M = M_1 \bigoplus M_2$, where $M_1 = \begin{matrix}
		2& & 3\\ &1&
	\end{matrix}$ and $M_2 = 1$. We have that $M$ is not $\tau$-rigid since $\tau(M) = M_2$, and this decomposition satisfies $M_1 \notin \Fac(M_2)$.
	
	 Let $\Gamma = \{(\mathsf{T}(M), M^{\bot}), (\mathsf{T}(1), 1^{\bot}) \}$ be the finest nested family such that $M_{\bullet} \in \mathcal{M}(\Gamma)$, where $M_{\bullet} = (M_j)_{j \in \mathbb{I}_2}$. Then $M_{\bullet} \in \mathcal{M}^{\dagger}(\Gamma)$. Since $M_2 \in \mathcal{T}_2$ and $M_1 \notin \mathcal{T}_2$, it follows that $\mathfrak{F}_1(M_{\bullet}) = f_2(M_1) = 2\bigoplus 3$. Theorem \ref{teoresma1} ensures that $(\{\Omega_1 = 2\bigoplus 3, \Omega_2 =1\}, \leq)$ is a pre-stratifying system, from which we obtain two stratifying systems induced by $M_{\bullet}$: $$(\{\Delta_1 = 3, \Delta_2 =1\}, \leq) \hbox{ and } (\{\Delta'_1 = 2, \Delta'_2 = 1\}, \leq).$$
\end{exem}

\section{A stratifying system of infinite size not indexed by $(\mathbb{N}, \leq)$}\label{sec7}

In \cite{sistemainfinito}, H. Treffinger demonstrated the existence of stratifying systems of infinite size indexed by $(\mathbb{N}, \leq)$ (where $\leq$ is the natural order on $\mathbb{N}$) in $n$-representation infinite algebras. The aim of this section is to show that, using the formalism of a nested family of torsion pairs, it is possible to obtain, all at once, infinite stratifying systems, including those constructed in \cite{sistemainfinito}. Moreover, we obtain stratifying systems of infinite size that cannot be indexed by $(\mathbb{N}, \leq)$, as shown below.

First, to facilitate readability, we will recall some key concepts. Let $n$ be a positive integer, and let $A$ be a connected, finite-dimensional $k$-algebra with global dimension at most $n$. Denote by $\mathcal{D}^b(\modd A)$ the bounded derived category of $\modd A$. The Nakayama functor will be denoted by $\nu$, and by $\nu_n := \nu \circ [-n]$, with $\nu_n^{-1}$ denoting its inverse. A finite-dimensional algebra with global dimension at most $n$ is called {\bf $\boldsymbol{n}$-representation infinite} if, for any indecomposable projective $A$-module $P$ we have that $\nu_n^{-j}(P) \in \modd A$ (i. e. concentrated in degree $0$) \cite{ninfinite}. 

An example of an $n$-representation infinite algebra is the Beilinson algebra $A = kQ/I$, where $Q$ is the following quiver:

\[\begin{tikzcd}
	1 & 2 & 3 & n & n+1
	\arrow["{a_0^1}", curve={height=-12pt}, from=1-1, to=1-2]
	\arrow["{a^1_n}"', curve={height=12pt}, from=1-1, to=1-2]
	\arrow["\vdots"{description}, shift left, draw=none, from=1-1, to=1-2]
	\arrow["{a_0^2}", curve={height=-12pt}, from=1-2, to=1-3]
	\arrow["{a^2_n}"', curve={height=12pt}, from=1-2, to=1-3]
	\arrow["\vdots"{description}, shift left, draw=none, from=1-2, to=1-3]
	\arrow["\cdots"{description}, draw=none, from=1-3, to=1-4]
	\arrow["{a^n_0}", curve={height=-12pt}, from=1-4, to=1-5]
	\arrow["{a_n^n}"', curve={height=12pt}, from=1-4, to=1-5]
	\arrow["\vdots"{description}, shift left, draw=none, from=1-4, to=1-5]
\end{tikzcd}\]
\noindent and $I$ is the ideal generated by $a_i^{k}a_j^{k+1}-a_j^{k}a_i^{k+1}$ for all $k \in \mathbb{I}_{n-1}$ and $i, j \in \mathbb{I}_n\cup\{0\}$ \cite[Example 2.15]{ninfinite}. Note that $\vert A_A \vert := N = n + 1$.

For the next theorem, we fix the totally ordered set $G = (\{-1, 1\} \times \mathbb{N}, \leqt)$, where $\leqt$ denotes the lexicographic order defined by $(y, m) \geet (x, l)$ if $y > x$, or $y = x$ and $m > l$.

\begin{teo}\label{x1x}
	Let $A$ be an $n$-representation infinite algebra with $n > 1$. Then, there exists a pre-stratifying system $\Omega$ of infinite size indexed by $G = (\{-1, 1\} \times \mathbb{N}, \leqt)$. Moreover, $\Omega$ cannot be indexed by $(\mathbb{N}, \leq)$.
\end{teo}
\begin{proof} Let $A$ be an $n$-representation infinite algebra such that $A_A = \bigoplus\limits_{j=1}^{N}P_j$, with $P_j$ indecomposable for all $j \in \mathbb{I}_{N}$. Define $M_{(x, l)} = \nu_n^{xl}(D^{\frac{x+1}{2}}(A))$ and $M_{\bullet} = (M_{(x, l)})_{(x,l) \in G}$, with the convention that $D^0(A) = A$.

For each $(y, m) \in G$, consider the torsion pair $(\mathcal{T}_{(y, m)}, \mathcal{F}_{(y, m)})$, where
$$
\mathcal{T}_{(y, m)} = \mathsf{T}\left(\{M_{(x, l)} \ \mid \ (x, l)\geqt (y, m) \}\right)$$ \text{and} $$\mathcal{F}_{(y, m)} = \left(\{M_{(x, l)} \ \mid \ (x, l)\geqt (y, m) \}\right)^{\bot}.$$
 Since $\Hom_A(M_{(y, m)}, M_{(x, l)}) = 0$ if $(y, m) \geet (x, l)$ by \cite[Proposition 2.3 (b), Proposition 4.10 (d)]{ninfinite}, Theorem \ref{galo} and Lemma \ref{lesma1} ensure that $\Gamma=\{(\mathcal{T}_{(y,m)}, \mathcal{F}_{(y,m)} )\}_{(y, m)\in G}$ is a nested family of torsion pairs, $M_{\bullet}\in\mathcal{M}^{*}(\Gamma)$ and $\mathfrak{F}(M_{\bullet}) = M_{\bullet}$.

Since $n > 1$, we have $\Ext_A^1(M_{(y, m)}, \mathfrak{F}_{(x, l)}(M_{\bullet})) = \Ext_A^1(M_{(y, m)}, M_{(x, l)}) = 0$ by \cite[Proposition 4.10 (f)]{ninfinite}, in particular if $(y, m) \geqt (x, l)$. Therefore, $M_{\bullet} \in \mathcal{M}^{\dagger}(\Gamma)$. According to Theorem~\ref{teoresma1}, $\Omega = \{M_{(y, m)}\}_{(y, m)\in G}$ is a pre-stratifying system. 

This pre-stratifying system cannot be indexed by $(\mathbb{N}, \leq)$. Indeed, there does not exist an order-preserving bijection $\sigma: G \longrightarrow (\mathbb{N},\leq)$; since $$\{(x, l) \leet (1, 1) \ \mid \ (x, l) \in G \}$$ is an infinite set, which, in turn, would imply that $\{n < \sigma(1, 1) \ \mid \ n \in \mathbb{N}\}$ is infinite, contradicting the well-ordering principle. 
\end{proof}

We can explicitly present all the stratifying systems constructed in the Theorem \ref{x1x}. To construct a stratifying system $\Delta$ from $\Omega$, we choose an {\edg indecomposable direct summand} $\Delta_{(x, l)}$ of $M_{(x, l)} = \bigoplus\limits_{j = 1}^{N} \nu_n^{xl}(D^{\frac{x+1}{2}}(P_j))$ for every $(x, l) \in G$. \cite[Proposition 4.10 (a), (b)]{ninfinite} ensures that $\nu_n^{xl}(D^{\frac{x+1}{2}}(P_j))$ is an indecomposable module for every $j \in \mathbb{I}_N$, which implies that $\Delta = \{\nu_n^{xl}(D^{\frac{x+1}{2}}(P_{f(x, l)}))\}_{(x, l) \in G}$ is a stratifying system for any function $f: G \longrightarrow \mathbb{I}_N$. In particular, the number of induced stratifying systems is in bijection with the number of functions $f: G \longrightarrow \mathbb{I}_N$.

\begin{cor} Let $A$ be an $n$-representation infinite algebra with $n > 1$ and $\vert A_A \vert > 1$. Then, there exist infinitely many stratifying systems of infinite size that cannot be indexed by $(\mathbb{N}, \leq)$. \end{cor} \qed

It is important to note that the order plays a central role in the definition of a stratifying system. Observe that $\Delta =(\{\Delta_1 = P(2), \Delta_2 = P(3)\}, \leq)$ and $\Delta' =(\{\Delta'_1 = P(3), \Delta'_2 = P(2)\})$ are distinct stratifying systems, despite defining the same set, where $P(2)$ and $P(3)$ are defined in Section \ref{x1x1}.

Since $G$ is a countable set, there exists a bijection $\sigma: \mathbb{N} \longrightarrow G$ (which does not preserve the order). We can define $N_{\bullet} = (N_k)_{k \in \mathbb{N}} \in \GR_{\mathbb{N}}(\modd A)$ by $N_k := M_{\sigma(k)}$. It is not possible to follow the proof of Theorem~\ref{x1x} to construct a stratifying system indexed by $(\mathbb{N}, \leq)$ induced by $N_{\bullet}$. In the proof, we heavily used the fact that $\Hom_A(\nu_n^{l}(DA), \nu_n^{-k}(A)) = 0$ for all $l, k \geq 0$ \cite[Proposition 4.10 (f)]{ninfinite}. On the other hand, there will exist infinitely many pairs of indices $(i, j) \in \mathbb{N}\times \mathbb{N}$ with $i > j$ such that $N_i = \nu_n^{-k}(A)$ and $N_j = \nu_n^{l}(DA)$, with $k, l \geq 0$. In general, we do not have that $\Hom_A(\nu_n^{-k}(A), \nu_n^{l}(DA)) = 0$.

\section{Example}\label{x1x1}

We conclude the paper with a detailed example illustrating Theorem~\ref{teoresma1}. In this example, we show how the stratifying system of size 5 presented in \cite[Remark 2.7]{stratifyingsimple} can be induced from an object $N_{\bullet} \in \mathcal{N}^{\dagger}(\Gamma)$, where $\Gamma$ is a suitable nested family of torsion pairs.

Let $A = kQ/I$, where $Q$ is the following quiver: \[\begin{tikzcd}
	3 & 1 & 2 & 4
	\arrow["\alpha", from=1-1, to=1-2]
	\arrow["\beta"', from=1-3, to=1-2]
	\arrow["\gamma"', from=1-4, to=1-3]
\end{tikzcd}\] and $I$ is the ideal generated by $\gamma \beta$. The Auslander-Reiten quiver of $A$ is given by:
\[\begin{tikzcd}
	& {P(3)} && {S(2)} && {S(4)} \\
	{P(1)} && {I(1)} && {P(4)} \\
	& {P(2)} && {S(3)}
	\arrow[from=1-2, to=2-3]
	\arrow[dotted, from=1-4, to=1-2]
	\arrow[from=1-4, to=2-5]
	\arrow[dotted, from=1-6, to=1-4]
	\arrow[from=2-1, to=1-2]
	\arrow[from=2-1, to=3-2]
	\arrow[from=2-3, to=1-4]
	\arrow[dotted, from=2-3, to=2-1]
	\arrow[from=2-3, to=3-4]
	\arrow[from=2-5, to=1-6]
	\arrow[from=3-2, to=2-3]
	\arrow[dotted, from=3-4, to=3-2]
\end{tikzcd}\]

Consider $N_{\bullet} = (N_k)_{k \in \mathbb{I}_5}$, where
$$N_1=P(1),\quad
N_2=P(2),\quad
N_3=I(1),\quad
N_4=P(4),\quad
N_5=S(4).$$
Note that $N := \bigoplus\limits_{k=1}^{5} N_k$ is not $\tau^{-}$-rigid, since there exists a non-zero morphism $N \rightarrow \tau^{-1}(N)$, where $\tau^{-1}(N) = I(1) \bigoplus S(3)$.

Consider the nested family of torsion pairs $\Gamma = \{(\mathcal{T}_k, \mathcal{F}_k)\}_{k \in \mathbb{I}_5}$, where
\begin{align*}
	\mathcal{T}_1 = \add\{P(2), P(3), P(4), I(1), S(2), S(3), S(4)\} \hspace{1.2cm} \mathcal{F}_1 = \add\{P(1)\},\\
	\mathcal{T}_2 = \add\{P(3), P(4), S(3), S(4)\} \hspace{1.95cm} \mathcal{F}_2 = \add\{P(1), P(2), S(2)\},\\
	\mathcal{T}_3 = \add\{P(4), S(3), S(4)\} \hspace{1.2cm} \mathcal{F}_3 = \add\{P(1), P(2), P(3), I(1), S(2)\},\\
	\mathcal{T}_4 = \add\{S(4)\} \hspace{1.2cm} \mathcal{F}_4 = \add\{P(1), P(2), P(3), P(4), I(1), S(2), S(3)\},\\
	\mathcal{T}_5 = \add\{0\} \hspace{0.88cm} \mathcal{F}_5 = \add\{P(1), P(3), P(2), P(4), I(1), S(2), S(3), S(4)\}.
\end{align*}

One checks that $N_{\bullet}\in\mathcal{N}(\Gamma)$. The substratum of $N_{\bullet}$ is given by $\mathfrak{T}(N_{\bullet}) = (t_{k-1}(N_k))_{k\in \mathbb{I}_5}$. Since the torsion pair $(\mathcal{T}_k, \mathcal{F}_k)$ is splitting if $k \neq 2$ we have $t_{k-1}(N_k) = N_k$ for $k \neq 3$. 

Note that there exists an exact sequence
\[0 \longrightarrow P(3) \longrightarrow I(1) \longrightarrow S(2) \longrightarrow 0,\]
where $P(3) \in \mathcal{T}_2$ and $S(2) \in \mathcal{F}_2$, hence $t_{2}(N_3) = P(3)$.

We also verify that $N_{\bullet}\in \mathcal{N}^{\dagger}(\Gamma)$. Indeed, we have that $$\Ext_A^1(\mathfrak{T}_k(N_{\bullet}), N_j) = \Ext_A^1(t_{k-1}(N_k), N_j) = 0$$ if $k \geq j$ since $t_{k-1}(N_k)$ is projective if $k \neq 5$, and the only non-split exact sequence ending at $S(4)$ is given by:
\[
0 \longrightarrow S(2) \longrightarrow P(4) \longrightarrow S(4) \longrightarrow 0.
\]
Since $\mathfrak{T}_k(N_{\bullet})$ is indecomposable for every $k \in \mathbb{I}_5$, Theorem~\ref{teoresma1} ensures that $\Delta = \{\Delta_k\}_{k \in \mathbb{I}_5}$ is a stratifying system, where $\Delta_k$ is $P(1)$, $P(2)$, $P(3)$, $P(4)$, and $S(4)$ for $k = 1, 2, 3, 4$, and $5$, respectively.

\begin{ack*}
	The authors are grateful to Eduardo Marcos and Hipolito Treffinger for their helpful suggestions.
	We also extend our thanks to the anonymous referees for their thorough review and constructive feedback,
	which significantly contributed to enhancing the quality of this work.
\end{ack*}

\begin{Funding*}
	The first named author acknowledges the financial support from DMAT-UFPR and CNPq Universal Grant 405540/2023-0. The second named author acknowledges the financial support from Capes. This study was financed in part by the Coordena\c{c}\~ao de Aperfei\c{c}oamento de Pessoal de N\'ivel Superior - Brasil (CAPES) - Finance Code  001.
\end{Funding*}

\bibliographystyle{mminimumeric}
\bibliography{bibliografia}

@article {MR4736629,
    AUTHOR = {Bautista, R. and P\'{e}rez, E. and Salmer\'{o}n, L.},
     TITLE = {Tame and wild theorem for the category of modules filtered by
              standard modules},
   JOURNAL = {J. Algebra},
  FJOURNAL = {Journal of Algebra},
    VOLUME = {650},
      YEAR = {2024},
     PAGES = {394--457},
      ISSN = {0021-8693},
   MRCLASS = {16G60 (16G20 16G70)},
  MRNUMBER = {4736629},
       DOI = {10.1016/j.jalgebra.2024.03.021},
       URL = {https://doi.org/10.1016/j.jalgebra.2024.03.021},
}

@article {MR1097029,
    AUTHOR = {Auslander, M. and Reiten, I.},
     TITLE = {Applications of contravariantly finite subcategories},
   JOURNAL = {Adv. Math.},
  FJOURNAL = {Advances in Mathematics},
    VOLUME = 86,
      YEAR = 1991,
    NUMBER = 1,
     PAGES = {111--152},
      ISSN = {0001-8708},
   MRCLASS = {16G10 (16D90)},
  MRNUMBER = 1097029,
MRREVIEWER = {Dieter Happel},
       DOI = {10.1016/0001-8708(91)90037-8},
       URL = {https://doi.org/10.1016/0001-8708(91)90037-8},
}

@article {MR1128706,
    AUTHOR = {C. Ringel},
     TITLE = {The category of modules with good filtrations over a
              quasi-hereditary algebra has almost split sequences},
   JOURNAL = {Math. Z.},
  FJOURNAL = {Mathematische Zeitschrift},
    VOLUME = 208,
      YEAR = 1991,
    NUMBER = 2,
     PAGES = {209--223},
      ISSN = {0025-5874},
   MRCLASS = {16G10 (16D90)},
  MRNUMBER = 1128706,
MRREVIEWER = {J. Antonio de la Pe\~{n}a},
       DOI = {10.1007/BF02571521},
       URL = {https://doi.org/10.1007/BF02571521},
}

@incollection {MR1211481,
    AUTHOR = {Dlab, V. and Ringel, C.},
     TITLE = {The module theoretical approach to quasi-hereditary algebras},
 BOOKTITLE = {Representations of algebras and related topics ({K}yoto,
              1990)},
    SERIES = {London Math. Soc. Lecture Note Ser.},
    VOLUME = 168,
     PAGES = {200--224},
 PUBLISHER = {Cambridge Univ. Press, Cambridge},
      YEAR = 1992,
   MRCLASS = {16E60 (16D90 16G20)},
  MRNUMBER = 1211481,
MRREVIEWER = {Otto Kerner},
}

@book{livroazul,
title =     {Elements of the representation theory of associative algebras Volume 1},
author =    {I. Assem and A. Skowronski and D. Simson},
publisher = {Cambridge University Press},
isbn =      {052158423X; 9780521584234; 9780521586313; 0521586313; 9780511345456; 0511345453; 9780521836104; 0521836107; 9780521544207; 0521544203; 9780521882187; 0521882184},
year =      {2006},
series =    {London Mathematical Society student texts 65, 71-72}}

@article{stratifying,
author={O. Mendoza and H. Treffinger},
title={Stratifying Systems Through $\tau$-Tilting Theory},
journal={Documenta Mathematica},
page={701-720},
year={2020},
volume={25}}

@article{stratifyingsimple,
title = {Stratifying systems via relative simple modules},
journal = {Journal of Algebra},
volume = {280},
number = {2},
pages = {472-487},
year = {2004},
issn = {0021-8693},
doi = {https://doi.org/10.1016/j.jalgebra.2004.06.018},
url = {https://www.sciencedirect.com/science/article/pii/S002186930400345X},
author = {E. Marcos and O. Mendoza and C. S\'aenz},
}

@article{exceptional,
title = {$\tau$-exceptional sequences},
journal = {Journal of Algebra},
volume = {585},
pages = {36-68},
year = {2021},
issn = {0021-8693},
doi = {https://doi.org/10.1016/j.jalgebra.2021.04.038},
url = {https://www.sciencedirect.com/science/article/pii/S0021869321002866},
author = {A. Buan and B. Marsh},
}

@article{reiten, 
title={$\tau$-tilting theory}, 
volume={150}, 
DOI={10.1112/S0010437X13007422}, 
number={3}, 
journal={Compositio Mathematica}, 
author={T. Adachi and O. Iyama and I. Reiten}, 
year={2014}, 
pages={415-452}
}

@book{coelho,
   title =     {Basic Representation Theory of Algebras (Graduate Texts in Mathematics (283), Band 283)},
   author =    {I. Assem and F. Coelho},
   publisher = {Springer},
   isbn =      {3030351173,9783030351175},
   year =      {2020},
   series =    {},
   edition =   {1},
   volume =    {}}

@article{sistemainfinito,
author = {H. Treffinger},
year = {2023},
month = {01},
pages = {15-19},
title = {The size of a stratifying system can be arbitrarily large},
volume = {361},
journal = {Comptes Rendus. Math\'ematique},
doi = {10.5802/crmath.385}
}

@phdthesis{tesecadavid,
    author = {P. Cadavid},
    title = {Sistemas estratificantes sobre \'algebras heredit\'arias},
    school = {Instituto de Matem\'atica e Estat\'istica, Universidade de S\~ao Paulo},
    link = {https://teses.usp.br/teses/disponiveis/45/45131/tde-05022015-111007/pt-br.php},
    year = 2012
}

@article{almostsplit,
title = {Almost split sequences in subcategories},
journal = {Journal of Algebra},
volume = {69},
number = {2},
pages = {426-454},
year = {1981},
issn = {0021-8693},
doi = {https://doi.org/10.1016/0021-8693(81)90214-3},
url = {https://www.sciencedirect.com/science/article/pii/0021869381902143},
author = {M. Auslander and S. Smal\o}
}

@article{ninfinite,
title = {n-representation infinite algebras},
journal = {Advances in Mathematics},
volume = {252},
pages = {292-342},
year = {2014},
issn = {0001-8708},
doi = {https://doi.org/10.1016/j.aim.2013.09.023},
url = {https://www.sciencedirect.com/science/article/pii/S0001870813003642},
author = {M. Herschend and O. Iyama and S. Oppermann},
}

@article{definicaosistemaestratificante,
author = {K. Erdmann and C. S\'aenz},
title = {On Standardly Stratified Algebras},
journal = {Communications in Algebra},
volume = {31},
number = {7},
pages = {3429--3446},
year = {2003},
publisher = {Taylor \& Francis},
doi = {10.1081/AGB-120022232},
URL = { https://doi.org/10.1081/AGB-120022232},
}

@article{Cline1988,
author = {E. Cline and B. Parshall and L. Scott},
journal = {Journal f$\ddot{u}$r die reine und angewandte Mathematik},
pages = {85-99},
title = {Finite dimensional algebras and highest weight categories.},
url = {http://eudml.org/doc/153071},
volume = {391},
year = {1988},
}

@Article{Quasi,
 Author = {V. Dlab and C. Ringel},
 Title = {Quasi-hereditary algebras},
 FJournal = {Illinois Journal of Mathematics},
 Journal = {Ill. J. Math.},
 ISSN = {0019-2082},
 Volume = {33},
 Number = {2},
 Pages = {280--291},
 Year = {1989},
 Language = {English},
 zbMATH = {4089744},
 Zbl = {0666.16014}
}

@article{Cadavid1,
author = {P. Cadavid and E. Marcos},
title = {Stratifying systems over hereditary algebras},
journal = {Journal of Algebra and Its Applications},
volume = {14},
number = {06},
pages = {1550093},
year = {2015},
doi = {10.1142/S0219498815500930},
}

@article{Cadavid2,
author = {P. Cadavid and E. Marcos},
title = {Stratifying systems over the hereditary path algebra with quiver $\mathbb {A}_{p,q}$},
journal = {S\~ao Paulo Journal of Mathematical Sciences},
volume = {10},
year = {2016},
doi = {10.1007/s40863-015-0029-x},
}

@article{jordan, 
title={Stratifying systems and Jordan-H\"{o}lder extriangulated categories}, 
DOI={10.1017/S0017089525100621}, 
journal={Glasgow Mathematical Journal}, 
author={T. Br\"ustle and S. Hassoun and A. Shah and A. Tattar}, 
year={2025}
}

@article{AGOSTON20084177,
	title = {Approximations of algebras by standardly stratified algebras},
	journal = {Journal of Algebra},
	volume = {319},
	number = {10},
	pages = {4177-4198},
	year = {2008},
	issn = {0021-8693},
	doi = {https://doi.org/10.1016/j.jalgebra.2008.02.017},
	url = {https://www.sciencedirect.com/science/article/pii/S0021869308001087},
	author = {I. \'Agoston and V. Dlab and E. Luk\'acs},
}

@article{Bruce,
	ISSN = {00029947},
	URL = {http://www.jstor.org/stable/1994341},
	author = {S. Dickson},
	journal = {Transactions of the American Mathematical Society},
	number = {1},
	pages = {223--235},
	publisher = {American Mathematical Society},
	title = {A Torsion Theory for Abelian Categories},
	urldate = {2024-12-24},
	volume = {121},
	year = {1966}
}

@article{HN,
	title = {An algebraic approach to Harder-Narasimhan filtrations},
	journal = {Journal of Pure and Applied Algebra},
	volume = {229},
	number = {1},
	pages = {107817},
	year = {2025},
	issn = {0022-4049},
	doi = {https://doi.org/10.1016/j.jpaa.2024.107817},
	url = {https://www.sciencedirect.com/science/article/pii/S0022404924002147},
	author = {H. Treffinger},

}









\end{document}